\documentclass[11pt]{article}
\usepackage{amssymb}
\usepackage{latexsym}
\usepackage{color}
\usepackage{amsmath}
\usepackage{amsfonts}
\usepackage{amsthm}
\usepackage{varioref}
\usepackage{pdfsync}
\usepackage{enumitem}
\usepackage{mathrsfs}
\usepackage[utf8]{inputenc}
\usepackage[T1]{fontenc}
\usepackage[english]{babel}
\DeclareUnicodeCharacter{014D}{\=o}
\usepackage{geometry}
\makeatletter \@addtoreset{equation}{section}

\def\bfX{{\mathbf X}}

\def\bfL#1{{\mathbf \Lambda^{#1}}}
\newcommand{\No}[1]{\left\|#1\right\|}     
\newcommand{\abs}[1]{\left|#1\right|}     

\def\bfx{\mathbf x}
\def\bfy{\mathbf y}

\def\Yw{\mathfrak Y}
\def\Zw{\mathfrak Z}

\def\abs#1{\left|#1\right|}
\def\norm#1{\left\|#1\right\|}

\def \Sum{\displaystyle\sum}

\def \E{\mathbb{E}}
\def \F{\mathbb{F}}

\def \N{\mathbb{N}}
\def \R{\mathbb{R}}

\def\P{\mathbb{P}}
\def\Q{\mathbb{Q}}

\def\Ac{{\cal A}}

\def\Ec{{\cal E}}
\def\Fc{{\cal F}}

\def\Lc{{\cal L}}
\def\Mc{{\cal M}}

\def\Pc{{\cal P}}

\def\Tc{{\cal T}}
\def\Uc{{\cal U}}
\def\Vc{{\cal V}}

\def\Yc{{\cal Y}}
\def\Zc{{\cal Z}}

\def\ep{\hbox{ }\hfill$\Box$}
\def\reff#1{{\rm(\ref{#1})}}
\def\be{\begin{eqnarray}}
\def\ee{\end{eqnarray}}
\def\bal{\begin{aligned}}
\def\eal{\end{aligned}}
\def\beq{\begin{equation}}
\def\eeq{\end{equation}}
\def\beqq{\begin{equation*}}
\def\eeqq{\end{equation*}}
\def\b*{\begin{eqnarray*}}
\def\e*{\end{eqnarray*}}

\def\x{\times}

\def\={\;=\;}

\def\pourtout{\mbox{ for all } }

\def\.{\;.}

\def\eps{\varepsilon}

\def\vp{\varphi}
\def\1{\mathbf1}

\def\Tr#1{\mbox{\rm Tr}\left[#1\right]}

\theoremstyle{plain}
\newtheorem{theorem}{Theorem}[section]
\newtheorem{proposition}[theorem]{Proposition}

\newtheorem{assumption}[theorem]{Assumption}
\newtheorem{definition}[theorem]{Definition}
\newtheorem{lemma}[theorem]{Lemma}
\theoremstyle{definition}
\newtheorem{remark}[theorem]{Remark}

\def\st{\text{ s.t. }}

\def \benumlab#1{\begin{enumerate}[label={\rm \bf{(#1{\arabic{*}})}}, ref={\rm (#1{\arabic{*}})}]}
\def \enumlab{\end{enumerate}}

\labelformat{assumptionCU}{(${\underline{{\rm\bf C}}}$)}
\labelformat{assumptionCV}{(${{{\rm\bf C}}}$)}

\def \benumlabi#1{\begin{enumerate}[label={\rm {(#1-\roman{*})}}, ref={\rm{(#1-\roman{*})}}]}
\def \enumlabi{\end{enumerate}}

\def\step#1#2{\noindent \emph{Step #1:~#2}}

\definecolor{darkred}{rgb}{0.8,0,0}
\definecolor{darkblue}{rgb}{0,0,0.7}
\definecolor{darkgreen}{rgb}{0,0.4,0}

\begin{document}

\title{On a class of path--dependent singular stochastic control problems
}
\author{Romuald Elie \thanks{Universit\'e Paris--Est Marne--la--Vall\'ee, romuald.elie@univ-mlv.fr. This author gratefully acknowledges the support of the ANR project Pacman, ANR-16-CE05-0027.}\and
     Ludovic Moreau
     \thanks{EY, ludovic.moreau@fr.ey.com.}
     \and Dylan Possama\"{i}
     \thanks{Universit\'e Paris--Dauphine, PSL Research University, CNRS, CEREMADE, 75016 Paris, France, possamai@cere-made.dauphine.fr. This author gratefully acknowledges the support of the ANR project Pacman, ANR-16-CE05-0027.
     }
         }

\date{\today}
\parindent=0pt

\maketitle \vspace{-1em}
\begin{abstract}
This paper studies a class of non--Markovian singular stochastic control problems, for which we provide a novel probabilistic representation. The solution of such control problem is proved to identify with the solution of a $Z-$constrained BSDE, with dynamics associated to a non singular underlying forward process. Due to the non-Markovian environment, our main argumentation relies on the use of comparison arguments for path dependent PDEs. Our representation allows in particular to quantify the regularity of the solution to the singular stochastic control problem in terms of the space and time initial data.  Our framework also extends to the consideration of degenerate diffusions, leading to the representation of the solution as the infimum of solutions to $Z-$constrained BSDEs. As an application, we study the utility maximisation problem with transaction costs for non--Markovian dynamics. 

\vspace{0.5em}
{\bf Keywords:} singular control, constrained BSDEs, path-dependent PDEs, viscosity solutions, transaction costs, regularity.
\end{abstract}

\section{Introduction}
 The study of singular stochastic problems initiated by Chernoff \cite{chernoff1968optimal} and  Bather and Chernoff \cite{bather1967sequential,bather1967sequential2} gave rise to a large literature, mostly motivated by its large scope of applications in economics or mathematics. This includes in particular the well--known monotone follower problem, see for instance Karatzas \cite{karatzas1982class}, real options decision modelling related to optimal investment issues, see Davis, Dempster, Sethi and Vermes \cite{davis1987optimal}, Dynkin games, see  Karatzas and Wang \cite{karatzas2001connections} or Boetius \cite{boetius2005bounded}, optimal stopping problems, see Karatzas and Shreve \cite{karatzas1984connections,karatzas1985connections}, Boetius and Kohlmann \cite{boetius1998connections}, or Benth and Reikvam \cite{benth2004connection}, as well as optimal switching problems, see Guo and Tomecek \cite{guo2008connections}. For all these questions, the problem of interest is naturally modelled under the form of a singular control of finite variation problem. In this abundant literature, it is quite striking that very few studies  take into consideration this type of questions in a non--Markovian framework. One of the purpose of his paper is to try to fill this gap\footnote{We learned while finishing this paper that Bouchard, Cheridito and Hu \cite{bouchard2017prep} were working on similar questions. Their work in preparation will deal with more general controlled SDEs, where the $f$ term in \eqref{eq: SDE weak} below will be allowed to depend on the state process. However, their method of proof will be purely probabilistic and thus of a completely different nature compared to ours. Furthermore, their approach will not, {\it a priori}, allow for degenerate diffusions}.

 \vspace{0.5em}
In a Markovian environment, the solution to nice singular stochastic control problems characterises typically as the unique weak solution to a variational inequality, where the linear part of the dynamics is combined with a constraint on the gradient of the solution. For example, when modelling the optimisation of dividend flow for a firm, as initiated in continuous time by Jeanblanc-Picqu\'e and Shiryaev \cite{jeanblanc1995optimization}, the optimal singular actions identify to paying dividends as soon as the underlying wealth process hits a free boundary, where the gradient of the value function reaches its upper--bound value $1$. This example naturally suggests a connection between singular stochastic control problems and stochastic processes with gradient constraints, and more precisely backward stochastic differential equations (BSDEs) with constraints on the gain process, as introduced by Cvitani\'c, Karatzas and Soner \cite{cvitanic1998backward}. The main initial motivation for the introduction of such type of equation was the super--hedging of claims with portfolio constraint. We establish in this paper that the solution to such BSDEs provides a nice probabilistic representation for solution to singular stochastic control problems.
 
  \vspace{0.5em}
More precisely, the class of non--Markovian stochastic control problem of interest is of the form 
\begin{equation*}
v^{\rm sing}: (t,\bfx) \longmapsto \underset{K\in\mathcal U^t_{\rm sing}}{\sup}\mathbb E^{\P_0^t}\big[U\big(\bfx\otimes_tX^{t,\bfx,K}\big)\big],\; 0\le t\le T,
\end{equation*} 
where $\mathcal U_{\rm sing}^t$ denotes the set of multidimensional c\`adl\`ag non decreasing $\F^{t}-$adapted process stating from $0$ and in $L^p$, for $p\ge1$. The controlled underlying process $X$ has the following non--Markovian singular dynamics 
    \begin{align*}
            X^{t,\bfx,K}=\bfx(t)+\int_t^{\cdot} \mu^{t,\bfx}_s\left(X^{t,\bfx,K}\right)ds+\int_t^{\cdot}f_sdK_s +\int_t^{\cdot} \sigma^{t,\bfx}_s\left(X^{t,\bfx,K}\right)dB^t_s,
    \end{align*}
 where $\mu$ and $\sigma$ are functional maps satisfying usual conditions, as detailed in Assumption \ref{ass: mu sigma lip} below. By density argument, it is worth noticing that we may also reduce to taking the supremum over the subset of absolutely continuous controls. We can also consider a weak version of such control problem.
 
  \vspace{0.5em}
 After a well suited Girsanov type probability transform, we rewrite $v^{\rm sing}(t,\bfx)$ in a form that looks similar to a a face--lift type transformation of the terminal reward. Hereby, this provides the intuition behind the representation of $v^{\rm sing}(t,\bfx)$ as the solution to a BSDE with $Z-$constraint. The convex constraints imposed on the integrand process $Z$ are induced by the directions $f$ and represented via the convex set 
        $$
            \mathfrak K_t:=\left\{ q\in\R^d \st (f_t^\top q)\cdot \mathbf e_i\le0 \pourtout i\in\{1,\dots,d\}\right\}\;.
        $$
For any time $t$, we verify that $v^{\rm sing}(t,\bfx)$ identifies to $\Yc^{t,\bfx}$ where $(\Yc^{t,\bfx},\Zc^{t,\bfx})$ is the minimal solution to the constrained BSDE
        \begin{align*}
            \Yc^{t,\bfx}_\cdot\ge U^{t,\bfx}\left( \bfX^{t,\bfx} \right)-\int_\cdot^T \Zc^{t,\bfx}_s\cdot dB^t_s,\;
            ((\sigma^{t,\bfx}_s)^\top)^{-1}\left(\bfX^{t,\bfx}\right)\Zc^{t,\bfx}_s\in K_s.
        \end{align*}
 The line of proof relies on the observation that penalised versions of both the singular problem and BSDE are solutions of the same path-dependent partial differential equation (PPDE for short), for which we are able to provide a comparison theorem. As far as we know, it is the first time that the newly introduced theory of viscosity solutions of path-dependent PDEs (see the works of Ekren, Keller, Ren, Touzi and Zhang \cite{ekren2014viscosity,ekren2016viscosity,ekren2012viscosity,ren2014comparison,ren2014overview,ren2015comparison}) is used to prove such a representation. Even though the ideas are reminiscent of the approach that one could have used in the Markovian case (see for instance Peng and Xu \cite{peng2007constrained} for probabilistic interpretation of classical variational inequalities through BSDEs with constraints), we believe that using it successfully in a non--Markovian setting will open the door to many new potential applications of the PPDE theory.
 
  \vspace{0.5em}
 Our main motivation for investigating such singular stochastic control problem, was the problem of utility maximisation faced by an investor in a market presenting both transaction costs and non--Markovian dynamics. From the point of view of applications, having non--Markovian dynamics can be seen as a very desirable effect, as this case encompasses stochastic volatility models for instance. However, in such a framework our previous representation does not apply directly, since it requires non-degeneracy of the diffusion matrix of $X$, and due to transaction costs, the dynamics of wealth needs to be viewed as a bi--dimensional stochastic process driven by a one dimensional noise. We therefore extend our representation in order to include degenerate volatility coefficients. Our line of proof relies on compactness properties together with convex order ordering arguments as in Pag\`es \cite{pages2014convex}. In such degenerate context, the solution to the singular control problem identifies with the infimum of a family of constrained BSDEs.
 
  \vspace{0.5em}
 As a by--product, the probabilistic representation of $v^{sing}$ in terms of a BSDE solution allows to derive insightful properties on the singular stochastic control problem. First, it automatically provides a dynamic programming principle for such problem. Second, this representation allows us to quantify the regularity of $v^{sing}$ in terms of the initial data points. We observe that $v^{sing}$ is Lipschitz in space as well as $1/2-$H\"older continuous in time. Obtaining such results for singular control problems is in general a very hard task (see the discussion in \cite[Section 4.2]{bouchard2016} for instance), and our approach could be one potential and promising solution.
 
  \vspace{0.5em}
The paper is organised as follows: Section \ref{Sec2} presents the class of singular control problems of interest and derives alternative mathematical representations. Section \ref{Sec3} presents the corresponding constrained BSDE representation. The connection is derived in Section \ref{Sec4} via path--dependent PDE arguments. The consideration of degenerate volatility together with the example on transaction costs example are discussed in Section \ref{Sec5}. Finally Section \ref{Sec6} provides the applications of the representation in terms of dynamic programming and regularity of the solution. 
 
 \vspace{0.5em}
{\bf Notations:} For any $d\in\mathbb N\backslash\{0\}$, and for every vector $x\in\mathbb R^d$, we will denote its entries by $x_i$, $1\leq i\leq d$. For any $p\geq l$ and $(x_l,x_{l+1},\dots,x_p)\in (\R^d)^{p-l+1}$, we will also sometimes use the notation $x_{l:p}:=(x_l,x_{l+1},\dots,x_p)$.

\section{The singular control problem}\label{Sec2}
\subsection{Preliminaries}

We fix throughout the paper a time horizon $T>0$. For any $(t,x)\in[0,T]\x\R^d$, we denote by $\bfL{t,x}$ the space of continuous functions $\bfx$ on $[0,T]$, satisfying $\bfx_t=x$, $B^{t,x}$ the corresponding canonical process and $\F^{t,x,o}:=(\Fc_s^{t,x,o})_{t\leq s\leq T}$ the (raw) natural filtration of $B^{t,x}$. It is classical result that the $\sigma$-algebra $\Fc_T^{t,x,o}$ coincides with the Borel $\sigma$-algebra on $\bfL{t,x}$, for the topology of uniform convergence. We will simplify notations when $x=0$ by setting $\Omega^t:=\bfL{t,0}$, $\Omega:=\Omega^0$, $B^{t}:=B^{t,0}$, $B:=B^0$, $\F^{t,o}:=\F^{t,0,o}$ and $\F^o:=\F^{0,o}$. Besides, we will denote generically by $\mathcal C^t$ the space of continuous functions on $[t,T]$, without any reference to their values at time $t$. Moreover, $\P^t_{x}$ will denote the Wiener measure on $(\bfL{t,x},\Fc_T^{t,x,o})$, that is the unique measure on this space which makes the canonical process $B^{t,x}$ a Brownian motion on $[t,T]$, starting from $x$ at time $t$. We will often make use of the completed natural filtration of $\F^{t,x,o}$ under the measure $\mathbb P^{t}_x$, which we denote $\F^{t,x}:=(\mathcal F_s^{t,x})_{t\leq s\leq T}$. Again we simplify notations by setting $\F^{t}:=\F^{t,0}$ and $\F:=\F^{0}$, and we emphasise that all these filtrations satisfy the usual assumptions of completeness and right-continuity. For any $t\in[0,T]$, any $s\in[t,T]$ and any $\bfx\in\mathcal C^t$, we will abuse notations and denote
$$\No{\bfx}_{\infty,s}:=\underset{t\leq u\leq s}{\sup}\No{\bfx_u},$$
where $\No{\cdot}$ is the usual Euclidean norm on $\mathbb R^d$, which we denote simply by $\abs{\cdot}$ when $d=1$. Furthermore, the usual inner product on $\mathbb R^d$ is denoted by $x\cdot y$, for any $(x,y)\in\mathbb R^d\times\mathbb R^d.$

\vspace{0.5em}
For any $(t,s)\in[0,T]\times[t,T]$ and any $\bfx\in\mathcal C^t$, we define $\bfx^s\in\mathcal C^s$ by
$$\bfx^s(r):=\bfx(r),\ r\in[s,T].$$

\vspace{0.5em}
We also define the following concatenation operation on continuous paths. For any $0\leq t< t'\leq s\leq T$, for any $(x,x')\in \mathbb R^d\times \mathbb R^d$ and any $(\bfx,\bfx')\in\bfL{t,x}\times\bfL{t',x'}$, we let $\bfx\otimes_s\bfx'\in\bfL{t,x}$ be defined as
$$(\bfx\otimes_s\bfx')(r):=\bfx(r)\mathbf{1}_{t\leq r\leq s}+\left(\bfx'(r)+\bfx(s)-\bfx'(s)\right)\mathbf{1}_{s<r\leq T}.$$
Let us consider some $(t,x,s,\bfx)\in[0,T]\times \mathbb R^d\times [t,T]\times\mathcal C^t$. We will also denote, for simplicity, by $x\otimes_s\bfx$ the concatenation between the constant path equal to $x$ on $[0,T]$ and $\bfx$. That being said, for any map $g:[0,T]\times\mathcal C^0$ and for any $(t,\bfx)\in[0,T]\times\mathcal C^0$, we will denote by $g^{t,\bfx}$ the map from $[t,T]\times\mathcal C^t$ defined by
$$g^{t,\bfx}(s,\bfx'):=g(s,\bfx\otimes_t\bfx').$$

Furthermore, we also use the following (pseudo)distance, defined for any $(t,t',s,s')\in[0,T]^4$, any $(x,x')\in \mathbb R^d\times \mathbb R^d$ and any $(\bfx,\bfx')\in\bfL{s,x}\times\bfL{s',x'}$ by
$$d_\infty\left((t,\bfx),(t',\bfx')\right):=\sqrt{\abs{t'-t}}+\underset{0\leq r\leq T}{\sup}\No{(x\otimes_s\bfx)(r\wedge t)-(x'\otimes_{s'}\bfx')(r\wedge t')}.$$

\subsection{A first version of the control problem}
The first set of control processes that we will consider will be typical of singular stochastic control. More precisely, we define
\begin{align*}
\mathcal U_{\rm sing}^t:=&\left\{(K_s)_{t\leq s\leq T}, \ \text{which are c\`adl\`ag, $\mathbb R^d$-valued, $\F^{t}$-adapted, null at $t$,}\right.\\
&\hspace{0.4em}\left. \text{ with non--decreasing entries and s.t. for any $p\geq 1$,}\  \mathbb E\left[K_T^p\right]<+\infty\right\}.
\end{align*}

We next consider the following maps
$$\mu:[0,T]\times \mathcal C^0\longmapsto \mathbb R^d\text{ and }\sigma:[0,T]\times\mathcal C^0\longmapsto \mathbb S^d,$$
where $\mathbb S^d$ is the set of $d\times d$ matrices (which we endow with the operator norm associated to $\No{\cdot}$, which we still denote $\No{\cdot}$ for simplicity) as well as a bounded map $f:[0,T]\longmapsto \mathbb S^d$.

\vspace{0.5em}
The following assumption will be in force throughout the paper
\begin{assumption}\label{ass: mu sigma lip}
\begin{itemize}
\item[$(i)$] The maps $\mu$ and $\sigma$ are progressively measurable, in the sense that for any $(\bfx,\bfx')\in\mathcal C^0\times\mathcal C^0$ and any $t\in[0,T]$, we have for $\varphi=\mu,\sigma$
$$\bfx(s)=\bfx'(s),\ \text{for all $s\in[0,t]$}\Rightarrow \varphi(s,\bfx)=\varphi(s,\bfx'), \ \text{for all $s\in[0,t]$}.$$
\item[$(ii)$] $\mu$ and $\sigma$ have linear growth in $\bfx$, uniformly in $t$, that is there exists a constant $C>0$ such that for every $(t,\bfx)\in[0,T]\times\mathcal C^0$
$$\No{\mu_t(\bfx)}+\No{\sigma_t(\bfx)}\leq C\left(1+\No{\bfx}_{\infty,t}\right).$$
\item[$(iii)$] $\mu$ and $\sigma$ are uniformly Lipschitz continuous in $\bfx$, that is there exists a constant $C>0$ such that for any $(t,\bfx,\bfx')\in[0,T]\times \mathcal C^0\times\mathcal C^0$ we have
$$\No{\mu_t(\bfx)-\mu_t(\bfx')}+\No{\sigma_t(\bfx)-\sigma_t(\bfx')}\leq C\No{\bfx-\bfx'}_{\infty,t}.$$
\item[$(iv)$] For any $(t,\bfx)\in[0,T]\times\mathcal C$, $\sigma_t(\bfx)$ is an invertible matrix and the matrix $\sigma_t^{-1}(\bfx)f$ is uniformly bounded in $(t,\bfx)$.
\end{itemize}
\end{assumption} 
For any $(t,\bfx)\in[0,T]\x\mathcal C^0$ and $K\in\Uc^t_{\rm sing}$, we denote respectively by $\bfX^{t,\bfx}$ and $X^{t,\bfx,K}$
    the unique strong solutions on $(\Omega^t,\mathcal F_T^{t,o},\mathbb P_0^t)$ of the following SDEs (existence and uniqueness under Assumption \ref{ass: mu sigma lip} are classical results which can be found for instance in \cite{jacod1979calcul}, see Theorems 14.18 and 14.21)
    \begin{align}
        \displaystyle
            \bfX^{t,\bfx}&=\bfx(t)+\int_t^{\cdot} \mu^{t,\bfx}_s\left(\bfX^{t,\bfx}\right)ds +\int_t^{\cdot} \sigma^{t,\bfx}_s\left(\bfX^{t,\bfx}\right)dB^t_s,
               \;\P_0^t-a.s.,\label{eq: SDE weak}\\
        \displaystyle
            X^{t,\bfx,K}&=\bfx(t)+\int_t^{\cdot} \mu^{t,\bfx}_s\left(X^{t,\bfx,K}\right)ds+\int_t^{\cdot}f_sdK_s +\int_t^{\cdot} \sigma^{t,\bfx}_s\left(X^{t,\bfx,K}\right)dB^t_s,
                \;\P_0^t-a.s.\nonumber
    \end{align}

 By Assumption \ref{ass: mu sigma lip},
                       and standard estimates on non-Markovian SDEs, we have the following lemma
                       \begin{lemma}\label{lem.estimates} For all $p\ge2$, there is $C_p>0$, depending only on $p$, $T$ and the constant $C$ in Assumption \ref{ass: mu sigma lip}, such that for all $(t,t';\bfx,\bfx';K,K')\in[0,T]\x[t,T]\x\bfL{}^2\x(\Uc^t_{\rm sing})^2$
                        \begin{align}
 &\displaystyle
                                \E^{\P^t_0}
                                \bigg[
                                    \sup_{t\le s\le t'}\norm{X^{t,\bfx,K}_s-\bfx(t)}^p
                                \bigg]\le C_p\Big((t'-t)^{\frac12}\big(1+\No{\bfx}^p_{\infty,t}\big)+\big(1+(t'-t)^\frac12\big)\E^{\P_0^t}\big[\No{K_{t'}-K_t}^{p}\big]\Big),
                                    \label{eq: estimate SDE Delta initial}
                            \\
                            &\displaystyle
                                \E^{\P^t_0}
                                \bigg[
                                    \sup_{t\le s\le T}\norm{X^{t,\bfx,K}_s}^p
                                \bigg]\le C_p\Big( 1+\No{\bfx}_{\infty,t}^p +\E^{\P_0^t}\big[\No{K_T-K_t}^{p}\big] \Big),
                                    \label{eq: estimate SDE growth}
                            \\
        &\displaystyle
                                \E^{\P^t_0}
                                \bigg[
                                    \sup_{t\le s\le T}\norm{X^{t,\bfx,K}_s-X^{t,\bfx',K'}_s}^p
                                \bigg]\le C_p\bigg(\No{\bfx-\bfx'}_{\infty,t}^p+\E^{\P_0^t}\bigg[\bigg\|\int_t^Td\big(K_s-K'_s\big)\bigg\|^{p}\bigg]\bigg).
                                    \label{eq: estimate SDE Delta}
                        \end{align}
\end{lemma}

The stochastic control problem we are interested in is then
\begin{equation}\label{eq:prob}
v^{\rm sing}(t,\bfx):=\underset{K\in\mathcal U^t_{\rm sing}}{\sup}\mathbb E^{\P_0^t}\big[U\big(\bfx\otimes_tX^{t,\bfx,K}\big)\big],
\end{equation}
where the reward function $U:\mathcal C^0\longrightarrow \R$ is assumed to satisfy
\begin{assumption}\label{assump:U}
For any $(\bfx,\bfx')\in\mathcal C^0\times\mathcal C^0$, we have for some $C>0$ and some $r\geq0$
$$\abs{U(\bfx)-U(\bfx')}\leq C\No{\bfx-\bfx'}_{\infty,T}\left(1+\No{\bfx}_{\infty,T}^r+\No{\bfx'}_{\infty,T}^r\right).$$
\end{assumption}

Notice that it is clear from \eqref{eq: estimate SDE growth} that under Assumption \ref{assump:U}, we have
$$\abs{v^{\rm sing}(t,\bfx)}<+\infty,\ \text{for any $(t,\bfx)\in[0,T]\times\mathcal C^0$.}$$

\subsection{Simplifications and reformulations of the problem}
The goal of this section is to first explain how we can without loss of generality restrict our attention to absolutely continuous control processes in the definition of $v^{\rm sing}$, and then to propose two useful reformulations of this simplified problem.

\vspace{0.5em}
Let us consider for any $t\in[0,T]$ the following subset $\mathcal U^t$ of $\mathcal U^t_{\rm sing}$ consisting of controls which are absolutely continuous with respect to the Lebesgue measure on $[t,T]$
\begin{align*}
\overline{\mathcal U}^t:=&\Big\{K\in\mathcal U^t_{\rm sing},\; K_s=\int_t^s\nu_rdr,\; \P_0^t-a.s.,\; \text{with $(\nu_s)_{t\leq s\leq T}$, $\F^t-$predictable and $(\R_+)^d-$valued}\Big\}.
\end{align*}
For any $K\in\overline{\mathcal U}^t$, it will be simpler for us to consider the corresponding process $\nu$, so that we define
\begin{align*}
\mathcal U^t:=&\bigg\{(\nu_s)_{t\leq s\leq T},\; \text{$(\R_+)^d-$valued, $\F^t-$predictable, and s.t. $\mathbb E^{\P_0^t}\bigg[\bigg\|\int_t^T\nu_sds\bigg\|^p\bigg]<+\infty$, $\forall p\geq 1$}\bigg\}.
\end{align*}
    Then, for any $(t,\bfx)\in[0,T]\x\mathcal C^0$ and $\nu\in\Uc^t$, we denote by $X^{t,\bfx,\nu}$
    the unique strong solution on $(\Omega^t,\mathcal F_T^{t,o},\mathbb P_0^t)$ of the following SDE
    \begin{align}\label{eq:sde}
            \displaystyle
            X^{t,\bfx,\nu}&=\bfx(t)+\int_t^{\cdot} \mu^{t,\bfx}_s\left(X^{t,\bfx,\nu}\right)ds+\int_t^{\cdot}f_s\nu_sds +\int_t^{\cdot} \sigma^{t,\bfx}_s\left(X^{t,\bfx,\nu}\right)dB^t_s,
                \;\P_0^t-a.s.
    \end{align}
    
We can then define
\begin{equation}\label{eq:prob2}
v(t,\bfx):=\underset{\nu\in\mathcal U^t}{\sup}\mathbb E^{\P_0^t}\big[U\big(\bfx\otimes_tX^{t,\bfx,\nu}\big)\big].
\end{equation}
Our next step is to introduce the so-called weak formulation for the control problem $v$. This basically boils down to to introduce, using Girsanov's theorem, probability measures equivalent to $\mathbb P_0^t$ and induced by control processes in $\mathcal U^t$, and to maximise over those measures. More precisely, for any $(t,\bfx)\in[0,T]\x\mathcal C^0$ and $\nu\in\Uc^t$, we now define the following $\P_0^t-$equivalent measure
    $$
        \frac{d\P^{t,\bfx,\nu}}{d\P_0^t}=\Ec\left( \int_t^\cdot (\sigma^{t,\bfx}_s)^{-1}\left( \bfX^{t,\bfx} \right)f_s\nu_s\cdot dB^t_s \right).
    $$
    The weak formulation of the control problem \eqref{eq:prob2} is defined as
    \begin{equation}\label{eq:prob3}
    v^{\rm weak}(t,\bfx):=\underset{\nu\in\mathcal U^t}{\sup}\; \mathbb E^{\mathbb P^{t,\bfx,\nu}}\big[U^{t,\bfx}\big(\bfX^{t,\bfx}\big)\big], \text{ for any $(t,\bfx)\in[0,T]\times\mathcal C^0$}.
    \end{equation}
Finally, our last ingredient will be to introduce what we coin a {\it canonical weak formulation} for $v$, which will be particularly well suited when we will use the theory of viscosity solutions for path--dependent PDEs. For any $(t,\bfx)\in[0,T]\times\mathcal C^0$, let us define the following probability measure on $(\Lambda^{t,\bfx_t},\Fc^{t,\bfx_t,o}_T)$
    $$\P_0^{t,\bfx}:=\P_0^t\circ\left(\bfX^{t,\bfx}\right)^{-1}.$$
    Since $\sigma$ is assumed to be invertible, it is a classical result that we have
    \begin{equation}\label{eq:filt}
    \F^{t,o}=\F^{\bfX^{t,\bfx}},\ \text{and }\F^t=\overline{\F^{\bfX^{t,\bfx}}}^{\mathbb P_0^t},
    \end{equation}
    where $\F^{\bfX^{t,\bfx}}$ denotes the raw natural filtration of $\bfX^{t,\bfx}$ and $\overline{\F^{\bfX^{t,\bfx}}}^{\mathbb P_0^t}$ its completion under $\mathbb P_0^t$. 

\vspace{0.5em}
Our main results of this section are summarised in the following proposition, whose proof is classical and postponed to the appendix.

\begin{proposition}\label{prop:simp}
Let Assumptions \ref{ass: mu sigma lip} and \ref{assump:U} hold. We have the following equalities for any $(t,\bfx)\in[0,T]\x\mathcal C^0$
$$v^{\rm sing}(t,\bfx)=v(t,\bfx)=v^{\rm weak}(t,\bfx)=\underset{\nu\in\mathcal U^t}{\sup}\;\E^{\P^{t,\bfx}_\nu}\left[U^{t,\bfx}\left(B^{t,\bfx_t}\right)\right].$$
\end{proposition}

\vspace{0.5em}
Given the above result, in the rest of the paper, we will therefore focus on the value function $v$ instead of $v^{\rm sing}$.

\begin{remark}
Notice that despite the fact that the result that maximisation under $\mathcal U_{\rm sing}^t$ and $\mathcal U^t$ actually lead to the same value function is classical, the continuity assumptions that we made on the functions intervening in our problem are crucial for it to hold. Indeed, as shown by Heinricher and Mizel \cite{heinricher1986stochastic}, such an approximation result is not always verified.
\end{remark}

\subsection{Approximating the value function}

To obtain our main probabilistic representation result for the value function $v$ (and thus for $v^{\rm sing}$), we will use, as mentioned before, the theory of viscosity solutions of path--dependent PDEs. However, in order to do so we will have to make a small detour, and first approximate $v$. 

\vspace{0.5em}
For any integer $n>0$ and any $t\in[0,T]$, we let $\mathcal U^{t,n}$ denote the subset of $\mathcal U^t$ consisting of processes $\nu$ such that $0\leq\nu_s^i\leq n$, for $i=1,\dots,d$, for Lebesgue almost every $s\in[t,T]$. We then define the approximating value function for all $(t,\bfx)\in[0,T]\times\mathcal C^0$
        \begin{align}\label{eq:vn_def}
            v^n(t,\bfx):=\sup_{\nu\in\Uc^{t,n}}\E^{\P^t_0}\big[ U^{t,\bfx}\big( X^{t,\bfx,\nu} \big) \big]
                    =\sup_{\nu\in\Uc^{t,n}}\E^{\P^{t,\bfx}_\nu}\big[U^{t,\bfx}(B^{t,\bfx})\big],
        \end{align}
 where the second equality can be proved exactly as for $v$ in Lemma \ref{prop:simp}. We have the following simple stability result, whose proof is postponed to the appendix.
 \begin{lemma}\label{lemma:conv}
 Under Assumptions \ref{ass: mu sigma lip} and \ref{assump:U}, for every $(t,\bfx)\in[0,T]\times\mathcal C^0$, we have that the sequence $(v^n(t,\bfx))_{n\geq 1}$ is non--decreasing and
 $$v^n(t,\bfx)\underset{n\rightarrow +\infty}{\longrightarrow}v(t,\bfx).$$
 \end{lemma}

    \section{The corresponding constrained BSDEs}\label{Sec3}
    \subsection{Spaces and norms}

        We now define the following family of convex sets, for any $t\in[0,T]$: 
        $$
            \mathfrak K_t:=\left\{ q\in\R^d \st (f_t^\top q)\cdot \mathbf e_i\le0 \pourtout i\in\{1,\dots,d\}\right\},
        $$
        where $({\mathbf e}_i)_{1\leq i\leq d}$ denotes the usual canonical basis of $\mathbb R^d$, and where for any $M\in\mathbb S^d$, $M^\top$ denotes its usual transposition. 
    \begin{remark} This form of constraint that one wishes intuitively to impose on the gradient of the value function $v$ is quite natural. Indeed, recall that $f$ describes the direction in which the underlying forward process is pushed in case of singular action.
    \end{remark}    
        
        We next introduce for any $p\geq 1$ the following spaces 
        \begin{align*}
      (i)\;  \mathbb S^p_{t}:=&\Big\{\text{$(Y_s)_{t\leq s\leq T}$, $\mathbb R-$valued, $\F^{t}-$progressively measurable, c\`adl\`ag, $\mathbb P_0^t-a.s.$,} \text{ s.t. $\No{Y}_{\mathbb S^p_t}<+\infty$}\Big\},
        \end{align*}
        where $$\No{Y}_{\mathbb S^p_t}^p:=\mathbb E^{\mathbb P^{t}_0}\bigg[\underset{t\leq s\leq T}{\sup}\abs{Y_s}^p\bigg].$$
                \begin{align*}
   (ii)\;     \mathfrak S^p_{t,\bfx}:=&\Big\{\text{$(Y_s)_{t\leq s\leq T}$, $\mathbb R-$valued, $\F^{t,\bfx}-$progressively measurable, c\`adl\`ag, $\mathbb P_0^{t,\bfx}-a.s.$,} \text{ $\No{Y}_{\mathfrak S^p_{t,\bfx}}<+\infty$}\Big\},
        \end{align*}
        where $$\No{Y}_{\mathfrak S^p_{t,\bfx}}^p:=\mathbb E^{\mathbb P^{t,\bfx}_0}\bigg[\underset{t\leq s\leq T}{\sup}\abs{Y_s}^p\bigg].$$
        
                \begin{align*}
     (iii)\;   \mathbb H^p_{t}:=&\Big\{\text{$(Z_s)_{t\leq s\leq T}$, $\mathbb R^d-$valued, $\F^{t}-$predictable, $\No{Z}_{\mathbb H^p_t}<+\infty$}\Big\},\; \No{Z}_{\mathbb H^p_t}^p:=\mathbb E^{\mathbb P^{t}_0}\Bigg[\bigg(\int_t^T\No{Z_s}^2ds\bigg)^{\frac p2}\Bigg].
        \end{align*}

                \begin{align*}
(iv)\;        \mathfrak H^p_{t,\bfx}:=&\Big\{\text{$(Z_s)_{t\leq s\leq T}$, $\mathbb R^d-$valued, $\F^{t,\bfx_t}-$predictable, $\No{Z}_{\mathfrak H^p_t}<+\infty$}\Big\},\; \No{Z}_{\mathfrak H^p_{t,\bfx}}^p:=\mathbb E^{\mathbb P^{t,\bfx}_0}\Bigg[\bigg(\int_t^T\No{Z_s}^2ds\bigg)^{\frac p2}\Bigg].
        \end{align*}

\subsection{Weak and strong formulations for the BSDE}
       For any $(t,\bfx)\in[0,T]\times\mathcal C^0$, we would like to solve the $\mathfrak K-$constrained BSDE with generator $0$ and terminal condition $U^{t,\bfx}\left(\bfX^{t,\bfx}\right)$, that is to say we want to find a pair $(\Yc^{t,\bfx},\Zc^{t,\bfx})\in \mathbb S^2_t\times\mathbb H^2_t$ such
        \begin{align}\label{eq:bsde1}
            \Yc^{t,\bfx}_\cdot\ge U^{t,\bfx}\left( \bfX^{t,\bfx} \right)-\int_\cdot^T \Zc^{t,\bfx}_s\cdot dB^t_s,\ \P^t_0-a.s.\\
            ((\sigma^{t,\bfx}_s)^\top)^{-1}\left(\bfX^{t,\bfx}\right)\Zc^{t,\bfx}_s\in \mathfrak K_s, \ ds\otimes d\P^t_0-a.e.,\label{eq:bsde2}
        \end{align}
        and such that if there is another pair $(\widetilde\Yc^{t,\bfx},\widetilde\Zc^{t,\bfx})\in \mathbb S^2_t\times\mathbb H^2_t$ satisfying \eqref{eq:bsde1} and \eqref{eq:bsde2}, then we have $\Yc^{t,\bfx}\leq \widetilde\Yc^{t,\bfx}$, $\mathbb P^t_0-a.s.$ When it exists, the pair $(\Yc^{t,\bfx},\Zc^{t,\bfx})$ is called the minimal solution of the $\mathfrak K-$constrained BSDE.
        
\vspace{0.5em}
It will also be important for us to look at the weak version of the constrained BSDE, where for any $(t,\bfx)\in[0,T]\times\mathcal C^0$, we now look for a pair $(\Yw^{t,\bfx},\Zw^{t,\bfx})\in \mathfrak S^2_{t,\bfx}\times\mathfrak H^2_{t,\bfx}$ such that
         \begin{align}\label{eq:bsde1w}
            \Yw^{t,\bfx}_\cdot\ge U^{t,\bfx}\left( B^{t,\bfx_t} \right)-\int_\cdot^T \Zw^{t,\bfx}_s\cdot dW^{t,\bfx}_s,\ \P^{t,\bfx}_0-a.s.\\
             ((\sigma^{t,\bfx}_s)^\top)^{-1}\left(B^{t,\bfx_t}\right)\Zw^{t,\bfx}_s\in \mathfrak K_s,\ ds\otimes d\P^{t,\bfx}_0-a.e.\label{eq:bsde2w},
        \end{align}
                and such that if there is another pair $(\widetilde\Yw^{t,\bfx},\widetilde\Zw^{t,\bfx})\in \mathfrak S^2_{t,\bfx}\times\mathfrak H^2_{t,\bfx}$ satisfying \eqref{eq:bsde1w} and \eqref{eq:bsde2w}, then we have $\Yw^{t,\bfx}\leq \widetilde\Yw^{t,\bfx}$, $\mathbb P^{t,\bfx}_0-a.s.$

\vspace{0.5em}
Before elaborating more on conditions ensuring existence of minimal solutions for the above equations, we will first spend some time explaining our main result and sketch its proof in the next section.

\subsection{Main result and sketch of proof}
Our main result explicits the links existing between the value function $v$ (and thus $v^{\rm sing}$) and the constrained BSDEs in weak and strong formulation introduced in the previous section.     \begin{theorem}\label{th:main}
Let Assumptions \ref{ass: mu sigma lip}, \ref{assump:U}, \ref{assump:exist} and \ref{assump:exist2} hold. Then, for any $(t,\bfx)\in[0,T]\times \mathcal C^0$, we have
$$v(t,\bfx)=\Yc^{t,\bfx}_t=\Yw_t^{t,\bfx}.$$
\end{theorem}        
        Our strategy of proof is inspired by an observation from the Markovian version of the control problem. It is of course known that the value function $v$ is in this case associated to a variational inequality with gradient constraint, which is nothing else than the corresponding Hamilton--Jacobi--Bellman equation (see for instance Fleming and Soner \cite[Chapter VIII]{fleming2006controlled}). Furthermore, the work of Peng and Xu \cite{peng2007constrained} has proved that Markovian $Z-$constrained BSDEs were naturally associated to viscosity solutions of specific types of variational inequalities with gradient constraint. Therefore, if these two objects both solve the same PDE, and the latter has a unique solution, then they must coincide.
        
        \vspace{0.5em}
        It is important to notice here that without looking at the PDE, the link between the singular control problem and the $Z-$constrained BSDE is far from being obvious. Indeed, even if it is known that BSDEs with $Z-$constraint have a solution which can be represented as the value function of a singular control problem (that is one of the crux of the seminal paper by Cvitani\'c and Karatzas \cite{cvitanic1998backward}), that control problem has a very specific form, involving in particular the so--called support function of the set where $Z$ is constrained (see Section \ref{sec:delta} below for details). There is therefore no reason in general for a given singular control problem like $v$ to be of this form. Indeed, even the set of admissible controls may actually be different for the two problems. The fact that bot control problems are linked to the same PDE is thus the cornerstone of the argument,  and it does not seem to be any immediate probabilistic argument available to obtain the same result. This is a situation reminiscent of the observation made by Bouchard \cite{bouchard2009stochastic2} that the PDEs appearing in the resolution of optimal switching problems as well as stochastic target problems with jumps were the same, and which inspired the link proved by Elie and Kharroubi between multidimensional obliquely--reflected BSDEs and one dimensional BSDEs with constrained jumps \cite{elie2014adding}, or the more recent representation, proved by Kharroubi and Pham \cite{kharroubi2015feynman}, of standard control problems in terms of BSDEs with constrained jumps. 
        
                \vspace{0.5em}
                Moving on to the non--Markovian case then seems hopeless, as the PDE argument described above is no longer available. This is exactly where the theory of PPDEs steps into play, as it provides us with an appropriate generalisation of the aforementioned Hamilton--Jacobi--Bellman equation. Nonetheless, the situation is not that simple, in the sense that a theory of path--dependent variational inequalities with gradient constraints, which would be our object of interest in a non--Markovian framework, does not exist yet. The reasons why are mainly of technical nature and stem from the fact that the associated control processes are unbounded, which destroys the crucial compactness properties needed (we refer the reader to the seminal papers on the subject \cite{ekren2014viscosity,ekren2016viscosity,ekren2012viscosity,ren2014comparison,ren2014overview,ren2015comparison} for more details). 
                
                                \vspace{0.5em}
                Therefore, our strategy of proof is as follows
                \begin{itemize}
                \item[$(i)$] First we approximate both the value function $v$ and the $Z-$constrained BSDE by their version with bounded controls, and ensure that the approximations do converge.
                       \item[$(ii)$] Prove using PPDE theory that both the approximations of the previous step are viscosity solutions of the same PPDE, for which we obtain a uniqueness result.
                              \item[$(iii)$] Pass to the limit and conclude.
                \end{itemize}
                We emphasise that even if the proof that the approximations converge are classical, and the proof that these approximations are viscosity solutions of the appropriate PPDEs is slightly more complicated than the one in the seminal papers netiowned above, our main contribution does not lie there. It lies in our approach tot he problem through the PPDE theory, which is completely novel, and on the fact that unlike the previous literature, which proved that $Z-$constrained BSDEs were associated to a variational inequality, which is also the Hamilton--Jacobi--Bellman equation of a specific singular control problem, we actually start with generic singular control problems, and prove that they are always linked to $Z-$constrained BSDEs.
        \subsection{Wellposedness for the constrained BSDEs}
The above constrained BSDEs have been studied in the literature, first in \cite{cvitanic1998backward} and then by Peng in \cite{peng1999monotonic}. However, all these existence results rely on the assumption that there is at least one solution (which does not have to be the minimal one) to the problem. This forces us to adopt the following assumptions.
\begin{assumption}\label{assump:exist}
For every $(t,\bfx)\in[0,T]\times \mathcal C^0$, there exists a pair $(\overline{\Yc}^{t,\bfx},\overline{\Zc}^{t,\bfx})\in \mathbb S^2_t\times\mathbb H^2_t$ such that  
\begin{align*}
            \overline{\Yc}^{t,\bfx}_\cdot\ge U^{t,\bfx}\left( \bfX^{t,\bfx} \right)-\int_\cdot^T  \overline{\Zc}^{t,\bfx}_s\cdot dB^t_s,\ \P^t_0-a.s.,\;  ((\sigma^{t,\bfx}_s)^\top)^{-1}\left(\bfX^{t,\bfx}\right) \overline{\Zc}^{t,\bfx}_s\in \mathfrak K_s, \ ds\otimes d\P^t_0-a.e.
                    \end{align*}
\end{assumption}
\begin{assumption}\label{assump:exist2}
For every $(t,\bfx)\in[0,T]\times \mathcal C^0$, there exists a pair $(\overline{\Yw}^{t,\bfx},\overline{\Zw}^{t,\bfx})\in \mathfrak S^2_{t,\bfx}\times\mathfrak H^2_{t,\bfx}$ such that  
\begin{align*}
            \overline{\Yw}^{t,\bfx}_\cdot\ge U^{t,\bfx}\left( B^{t,\bfx_t} \right)-\int_\cdot^T  \overline{\Zw}^{t,\bfx}_s\cdot dW^{t,\bfx}_s,\ \P^{t,\bfx}_0-a.s.,\; ((\sigma^{t,\bfx}_s)^\top)^{-1}\left(B^{t,\bfx_t}\right) \overline{\Zw}^{t,\bfx}_s\in \mathfrak K_s, \ ds\otimes d\P^{t,\bfx}_0-a.e.
                    \end{align*}
\end{assumption}

\begin{remark}These assumptions simply indicates that it is indeed possible to satisfy the $Z-$constraint as well as solve the BSDE. Such constrained BSDEs have been first introduced in order to find the super--hedging price of a claim under portfolio constraints, and such condition in this framework simply indicates that one can find an admissible portfolio strategy that indeed super--hedges the claim of interest. 
\end{remark}
We then have immediately from \cite{cvitanic1998backward} the following
\begin{proposition}\label{prop:exist}
Let Assumptions \ref{ass: mu sigma lip}, \ref{assump:U}, \ref{assump:exist}, and \ref{assump:exist2} hold. Then, the minimal solution of the $\mathfrak K-$constrained BSDEs \eqref{eq:bsde1} and \eqref{eq:bsde2} exist and we have that
$$\text{the law of $\Yc^{t,\bfx}$ under $\mathbb P_0^t$ $=$ the law of $\Yw^{t,\bfx}$ under $\mathbb P^{t,\bfx}_0$}.$$
\end{proposition}

\proof
We only argue for the strong formulation, as the the arguments are exactly similar for the weak formulation. Since it is clear by Assumption \ref{assump:U} and \eqref{eq: estimate SDE growth} that $\mathbb E^{\P^t_0}\big[\big | U^{t,\bfx}\big(\bfX^{t,\bfx}\big)\big|^2\big]<+\infty$, the result is an immediate consequence of the main result in \cite{cvitanic1998backward}. The equality in law is then clear by definition.
\ep

\vspace{0.5em}
\begin{remark}
Since Assumptions \ref{assump:exist} and \ref{assump:exist2} are rather implicit, let us discuss some sufficient conditions under which they hold. We can for instance use Assumption 7.1 in \cite{cvitanic1998backward}, which states that if there exist a constant $C\in\R$ and a process $\varphi\in\mathbb H^2_t$ such that $ ((\sigma^{t,\bfx}_s)^\top)^{-1}\left(\bfX^{t,\bfx}\right) \varphi_s\in \mathfrak K_s, \ ds\otimes d\P^t_0-a.e.$ and such that
\begin{equation}\label{eq:domination}
U^{t,\bfx}(\bfX^{t,\bfx})\leq C+\int_t^T\varphi_s\cdot dB^t_s,\ \mathbb P^t_0-a.e.,\end{equation}
then Assumption \ref{assump:exist} is satisfied. A similar condition can of course be written for the weak formulation. Notice also that \eqref{eq:domination} holds if for instance $U$ is bounded.

\vspace{0.5em}
Finally, notice that as an immediate consequence of the Blumenthal $0-1$ law and Proposition \ref{prop:exist}, we have the following equality for every $(t,\bfx)\in[0,T]\times\mathcal C^0$
$$\Yc^{t,\bfx}_t=\Yw^{t,\bfx}_t.$$
\end{remark}

   \subsection{The penalised BSDEs}

Exactly as we have approximated the value function $v$ by $v^n$ defined in \eqref{eq:vn_def}, it will be useful for us to consider approximations of the $\mathfrak K-$constrained BSDEs introduced in the previous section. It is actually a very well--known problem, which already appeared in \cite{cvitanic1998backward}, and which can be solved by considering the so--called penalised BSDEs associated to the $\mathfrak K-$constrained BSDE. Before doing so, we need to introduce, for any $(t,\bfx)\in[0,T]\x\mathcal C^0$, the map 
\beq\label{eq definition rho}
        \rho:q\in\R^d\longmapsto q^+\cdot\mathbf1_d
    \eeq
    where for each $q:=(q_1,\cdots,q_d)\top\in\R^d$ we have used the notation: $q^+:=(q_1^+,\cdots,q_d^+)^\top$.

\vspace{0.5em}
Under Assumptions \ref{ass: mu sigma lip} and \ref{assump:U}, we can then define for any $(t,\bfx,n)\in[0,T]\x\mathcal C^0\x\N^\ast$, $(\Yc^{t,\bfx,n},\Zc^{t,\bfx,n})\in\mathbb S^2_{t}\times\mathbb H^2_t$ as the unique solution of the following BSDE
        \beq\label{eq: penalized BSDE}
            \Yc^{t,\bfx,n}_\cdot=
                U^{t,\bfx}\left( \bfX^{t,\bfx} \right)
                +\int_\cdot^T n\rho\big( f_s^\top\big(\big(\sigma^{t,\bfx}_s\big)^\top\big)^{-1}\big(\bfX^{t,\bfx}\big)\Zc^{t,\bfx,n}_s \big)ds
                -\int_\cdot^T \Zc^{t,\bfx,n}_s dB^t_s,\
               \P^t_0-a.s.
        \eeq
        
        Notice that existence and uniqueness hold using for instance the results in \cite{el1997backward}, since under Assumptions \ref{ass: mu sigma lip} and \ref{assump:U}, the terminal condition is obviously square--integrable, the generator $z\longmapsto\rho(f_s^\top((\sigma^{t,\bfx})^\top)^{-1}(\bfX^{t,\bfx})z )$ is null at $0$ and uniformly Lipschitz continuous in $z$ (we remind the reader that $\sigma^{-1}f$ is bounded, so its transpose is bounded as well).
        
        \vspace{0.5em}
         Alternatively, we may consider the penalised BSDEs in weak formulation
        $$
            \Yw^{t,\bfx,n}_\cdot= U^{t,\bfx}\left( B^{t,\bfx} \right)
                +\int_\cdot^T n\rho\big( f_s^\top((\sigma^{t,\bfx}_s)^\top)^{-1}\left(B^{t,\bfx}\right)\Zw^{t,\bfx,n}_s \big)ds
                -\int_\cdot^T \Zw^{t,\bfx,n}_s dW^{t,\bfx}_s,
                \; \P^{t,\bfx}_0-a.s.,
        $$
        which also admit a unique solution under Assumptions \ref{ass: mu sigma lip} and \ref{assump:U}. We then have the following classical result (see \cite{cvitanic1998backward} or \cite{peng1999monotonic} for the proof).

\begin{lemma}\label{lemma:limit}
Let Assumptions \ref{ass: mu sigma lip}, \ref{assump:U}, \ref{assump:exist}, and \ref{assump:exist2} hold. Then, for any $(t,\bfx)\in[0,T]\times\mathcal C^0$
        \begin{align*}
     &   \Yc^{t,\bfx,n}_s\underset{n\rightarrow+\infty}{\uparrow}\Yc^{t,\bfx}_s,\ \text{for any $s\in[t,T]$},\ \mathbb P^t_0-a.s.,\ \text{and }\No{\Yc^{t,\bfx}-\Yc^{t,\bfx,n}}_{\mathbb H^2_t}\underset{n\rightarrow+\infty}{\longrightarrow}0,\\
     &   \Yw^{t,\bfx,n}_s\underset{n\rightarrow+\infty}{\uparrow}\Yw^{t,\bfx}_s,\ \text{for any $s\in[t,T]$},\ \mathbb P^{t,\bfx}_0-a.s.,\ \text{and }\No{\Yw^{t,\bfx}-\Yw^{t,\bfx,n}}_{\mathfrak H^2_{t,\bfx}}\underset{n\rightarrow+\infty}{\longrightarrow}0.
        \end{align*}       
\end{lemma}

\section{Proof of Theorem \ref{th:main}}\label{Sec4}
The goal of this section is to prove the representation
$$v(t,\bfx)=\Yc_t^{t,\bfx}=\Yw_t^{t,\bfx},\ \text{for every $(t,\bfx)\in[0,T]\times\mathcal C^0$.}$$

In order to prove this result, we will show that both $v^n(t,\bfx)$ defined in \eqref{eq:vn_def} and $u^n(t,\bfx):=\Yw^{t,\bfx,n}_t$ are viscosity solutions of a semi-linear path-dependent PDE, for which a comparison result holds. It will then imply that $v^n(t,\bfx)=u^n(t,\bfx)$, and the desired result will be obtained by passing to the limit when $n$ goes to $+\infty$, see Lemmas \ref{lemma:conv} and \ref{lemma:limit}.

\subsection{A crash course on PPDEs}

In this section, we follow closely \cite{ren2014comparison} to introduce all the notions needed for the definition of viscosity solutions of path-dependent PDEs. Let us start with the notions of regularity we will consider.

\begin{definition}
$(i)$ For any $(t,\bfx,\widetilde\bfx)\in[0,T]\times\mathcal C^0\times\mathcal C^t$, any $s\in[t,T]$ and any $d\geq 1$, we say that an $\mathbb R^d-$valued process $u$ on $(\bfL{}^{t,\bfx_t},\mathcal F_T^{t,\bfx_t,o})$ is in $C^0\big([s,T]\times\bfL{}^{s,(\bfx\otimes_t\widetilde\bfx)_s},\mathbb R^d\big)$ when it is continuous with respect to the distance $d_\infty$, that is to say that for any $\varepsilon>0$, for any $(r_1,r_2,\bfx_1,\bfx_2)\in[s,T]^2\times\mathcal C^s\times\mathcal C^s$, there exists $\delta>0$ such that 
$$\text{if $d_\infty((r_1,\bfx_1),(r_2,\bfx_2))\leq \delta \; \Longrightarrow \; \big\|u^{s,\widetilde\bfx}(r_1,\bfx_1)-u^{s,\widetilde\bfx}(r_2,\bfx_2)\big\|\leq \varepsilon$.}$$

\vspace{0.5em}
$(ii)$ For any $(t,\bfx,\widetilde\bfx)\in[0,T]\times\mathcal C^0\times\mathcal C^t$ and any $s\in[t,T]$, we say that an $\mathbb R-$valued process $u$ on $(\bfL{}^{t,\bfx_t},\mathcal F_T^{t,\bfx_t,o})$ belongs to $C^{1,2}\big([s,T]\times\bfL{}^{s,(\bfx\otimes_t\widetilde\bfx)_s}\big)$ if $u\in C^0\big([s,T]\times\bfL{}^{s,(\bfx\otimes_t\widetilde\bfx)_s},\R\big)$ and if there exists  $(Z,\Gamma)\in C^0\big([t,T]\times\bfL{}^{t,\bfx_t},\mathbb R^d\big)\times C^0\big([t,T]\times\bfL{}^{t,\bfx_t},\mathbb R\big)$ such that
$$u^{s,\widetilde\bfx}_z-u_s(\widetilde\bfx)=\int_s^z\Gamma^{s,\widetilde\bfx}_rdr+\int_s^zZ_r^{s,\widetilde\bfx}\cdot dB^{s,(\bfx\otimes_t\widetilde\bfx)_s}_r,\ z\in[s,T],\ \mathbb P_0^{s,(\bfx\otimes_t\widetilde\bfx)}-a.s.$$
We then denote for simplicity and for any $\widetilde \bfx\in\mathcal C^t$ and any $s\in[t,T]$, 
$$
\mathcal L^{t,\bfx}u(s,\widetilde\bfx):=\Gamma_s(\widetilde\bfx),\; Du(s,\widetilde\bfx):=((\sigma_s^{t,\bfx}(\widetilde\bfx))^\top)^{-1}Z_s(\widetilde\bfx).$$
\end{definition}

Let us then denote, for any $(t,\bfx)\in[0,T]\times\mathcal C^0$, by $\mathcal T^{t,\bfx}$ the set of $\F^{t,\bfx_t,o}-$stopping times taking values in $[t,T]$, by $\mathcal T^{t,\bfx}_+$ the subset of $\mathcal T^{t,\bfx}$ consisting of the stopping times taking values in $(t,T]$, and for any $H\in\mathcal T^{t,\bfx}$, by $\mathcal T^{t,\bfx}_H$ and $\mathcal T^{t,\bfx}_{H,+}$, the subsets of $\mathcal T^{t,\bfx}$ consisting of stopping times taking values respectively in $[t,H]$ and $(t,H]$.

\vspace{0.5em}
    Next, we define for any $N\geq 1$
    $$
        \Pc^{t,\bfx,N}:=\left\{ \P^{t,x}_\nu,\ \nu\in\Uc^{t,N} \right\},$$
        $$\Mc^{t,\bfx,N}:=\left\{ \Q\; s.t.\; \frac{d\Q}{d\P_0^{t,\bfx}}=\Ec\left(\int_t^Tb_sdW_s^{t,\bfx}\right),\; \text{$b$, $\F^{t,\bfx_t}-$predictable s.t. $\norm{b}_\infty\leq N$} \right\}.
    $$
    For any $(t,\bfx)\in[0,T]\x\mathcal C^0$, and for any $w\in C^0\big([0,T]\times\bfL{}^{0,\bfx_0},\R\big)$,
    we now define the sets of test functions for $w$ as
    \begin{align*}
                &\displaystyle \overline\Ac^n w(t,\bfx):=\left\{
                        \displaystyle \vp\in C^{1,2}([t,T]\times\bfL{t,\bfx_t}),\ 0=\left( \vp-w^{t,\bfx} \right)(t,\bfx^t)>
                        \displaystyle
                                \overline\Ec^n_{t}\left[ \left( \vp-w^{t,\bfx} \right)
                                    \left(
                                       \cdot,B^{t,\bfx_t}
                                    \right)_{\tau\wedge H}
                                    \right]\right.\\
                        &\left.\displaystyle \hspace{6.2em} \mbox{for some } H\in \mathcal T^{t,\bfx} \text{ and for all }\tau\in\mathcal T^{t,\bfx}_{H,+}
                        \right\},\\[0.8em]
                 &\displaystyle \underline\Ac^n w(t,\bfx)
                    :=\left\{
                        \displaystyle \vp\in C^{1,2}([t,T]\times\bfL{t,\bfx_t}),\ 0=\left( \vp-w^{t,\bfx} \right)(t,\bfx^t)<
                        \displaystyle
                                \underline\Ec^n_{t}\left[ \left( \vp-w^{t,\bfx} \right)
                                    \left(
                                       \cdot,B^{t,\bfx_t}
                                    \right)_{\tau\wedge H}
                                    \right]\right.\\
                        &\left.\displaystyle \hspace{6.2em} \mbox{for some } H\in \mathcal T^{t,\bfx} \text{ and for all }\tau\in\mathcal T^{t,\bfx}_{H,+}
                        \right\},
    \end{align*}
  where for all $\mathcal F^{t,\bfx_t}_T-$measurable $\xi$ such that the quantities below are finite
    $$
                \overline\Ec^n_t[\xi]:=\sup_{\Q\in\Mc^{t,\bfx,\bar n}}\E^{\Q}[\xi]\;, \;
                \underline\Ec^n_t[\xi]:=\inf_{\Q\in\Mc^{t,\bfx,\bar n}}\E^{\Q}[\xi]\;,
                \mbox{ with } \bar n:=n\sqrt{d}\max\left(1;\norm{\sigma^{-1}f}\right)\;.
    $$
    Finally, we define for every $(t,\bfx,\vp)\in[0,T]\x\mathcal C^0\x C^{1,2}\big([t,T]\times\bfL{}^{t,\bfx_t}\big)$ the following PPDE
    \beq\label{eq: penalized PPDE}
        -\Lc^{t,\bfx}\vp(t,\bfx^t)-n\rho\Big(f_t^\top D\vp(t,\bfx^t)\Big)=0.
    \eeq
    
    \begin{definition}
  Fix some $x\in\mathbb R^d$ and let $u\in C^0\big([0,T]\times\bfL{}^{0,x},\mathbb R\big)$. We say that
  
  \vspace{0.5em}
  $(i)$ $u$ is a viscosity sub--solution of {\rm PPDE} \eqref{eq: penalized PPDE} if for any $(t,\bfx,\varphi)\in[0,T)\times\mathcal C^0\times \underline\Ac^n u(t,\bfx)$
  $$  -\Lc^{t,\bfx}\vp(t,\bfx^t)-n\rho\Big(f_t^\top D\vp(t,\bfx^t)\Big)\leq0.$$
  
    $(ii)$ $u$ is a viscosity super--solution of {\rm PPDE} \eqref{eq: penalized PPDE} if for any $(t,\bfx,\varphi)\in[0,T)\times\mathcal C^0\times \overline\Ac^n u(t,\bfx)$
  $$  -\Lc^{t,\bfx}\vp(t,\bfx^t)-n\rho\Big(f_t^\top D\vp(t,\bfx^t)\Big)\geq0.$$
  $(iii)$ $u$ is a viscosity solution of {\rm PPDE} \eqref{eq: penalized PPDE} if it is both a sub-- and a super--solution.
    \end{definition}
    We shall end this section with the following result which will be useful in the PPDE derivation of the value function. The technical proof is postponed to the appendix.
    \begin{lemma}\label{lemma:utile}
      For all $(t,\bfx)\in[0,T]\x\mathcal C^0$ and $\tau\in\Tc^{t,\bfx}_+$ we have
      $$
        \underline{\Ec}^n_t[\tau-t]>0.
      $$
    \end{lemma}
   
    \subsection{The main result}
    
We first start by linking the approximated value function and BSDEs to PPDE \ref{eq: penalized PPDE}. Notice that even though the arguments used are classical, our estimates are substantially more involved as we assumed only local Lipschitz--continuity of our coefficients.

 \begin{proposition}\label{prop: ppde for app vf}
Under Assumptions \ref{ass: mu sigma lip}, \ref{assump:U}, \ref{assump:exist} and \ref{assump:exist2}, $v^{n}$ and $u^{n}(t,\bfx):=\Yw^{t,\bfx,n}_t=\Yc^{t,\bfx,n}_t$ are viscosity solutions of {\rm PPDE} \reff{eq: penalized PPDE}.
        \end{proposition}

Define now, for any $\bfx\in\mathcal C^0$, the following subset of $C^0\big([0,T]\times\bfL{}^{0,\bfx_0}\big)$
\begin{align*}
C^0_2\big([0,T]\times\bfL{}^{0,\bfx_0}\big):=&\Big\{u\in C^0\big([0,T]\times\bfL{}^{0,\bfx_0}\big),\ \text{s.t. for any $(t,\widetilde\bfx)\in[0,T]\times\mathcal C^0$,}\\
&\hspace{1em} \text{$u^{t,\bfx}$ is continuous in time $\mathbb P_0^{t,\widetilde\bfx}-a.s.$},\ u^{t,\widetilde\bfx}\in\mathfrak S^2_{t,\widetilde\bfx} \Big\}.
\end{align*}
We now recall the following comparison theorem from \cite{ren2014comparison} (see their Theorem 4.1), adapted to our context.
\begin{theorem}[\cite{ren2014comparison}]\label{th:comp}
Let $u,v$ in $C^0_2\big([0,T]\times\bfL{}^{0,\bfx_0}\big)$ be respectively viscosity sub--solution and super--solution of {\rm PPDE} \eqref{eq: penalized PPDE}. If $u(T,\cdot)\leq v(T,\cdot)$, then $u\leq v$ on $[0,T]\times\mathcal C^0$.
\end{theorem}

From Proposition \ref{prop: ppde for app vf}, we know that for every $n\geq 1$, $v^n$ and $u^n$ are viscosity solutions of PPDE \eqref{eq: penalized PPDE}. Since it is clear by all our estimates that $v^n,u^n\in C^0_2\big([0,T]\times\bfL{}^{0,\bfx_0}\big)$, and since $v^n(T,\cdot)=u^n(T,\cdot)$, by Theorem \ref{th:comp} we deduce that
$$v^n(t,\bfx)=\Yc^{t,\bfx,n}_t=\Yw_t^{t,\bfx,n}.$$
By Lemmas \ref{lemma:conv} and \ref{lemma:limit}, it then suffices to let $n$ go to $+\infty$ to conclude the proof of Theorem \ref{th:main}.
    \section{Extension to degenerate diffusions}\label{Sec5}
    \subsection{The setting}
    The result of the previous section is fundamentally based on the non--degeneracy of the diffusion matrix $\sigma$. Our main purpose here is to extend our general representation to cases where $\sigma$ is allowed to be degenerate. As will be clear later on, the type of degeneracy we will consider will be rather specific, but it will nonetheless be particularly well--suited for the applications we have in mind. Before stating our results, we need to introduce some notations. For every $n\in\mathbb N\backslash\{0\}$ and any $t\in[0,T)$, we consider uniform partitions of the interval $[t,T]$ by $\{t^{t,n}_k:=t+k(T-t)n^{-1},\ k=0,\dots,n\}$. We also define for every $0\leq k\leq n$ and every $(s,x_{0:k})\in[t,T]\times(\R^d)^{k+1}$, the linear interpolator $i_k:(\R^d)^{k+1}\longrightarrow \mathcal C^0$ by
  $$i_k(x_{0:k})(s)= \frac{n}{T-t}\sum_{i=0}^{k-1}\left((t^{t,n}_{i+1}-s)x_i+(s-t^{t,n}_i)x_{i+1}\right){\bf 1}_{[t^{t,n}_i,t^{t,n}_{i+1}]}(s).$$
Our main assumption now becomes
    \begin{assumption}\label{ass: mu sigma lip2}
    Assumption \ref{ass: mu sigma lip}$(i),(ii),(iii)$ hold and

\vspace{0.5em}
$(iv')$ For any $p> 0$, there exist progressively measurable maps $\eta^p:[0,T]\times\mathcal C^0\longrightarrow \R_-$ with linear growth, that is there exists some $C>0$ such that
$$0\leq-\eta^p_t(\bfx)\leq C(1+\No{\bfx}_{\infty,t}),$$
 and a deterministic map $M:[0,T]\longrightarrow \mathbb S^d$, such that $M_t$ is symmetric positive for every $t\in[0,T]$, the maps $\bfx\longmapsto\eta^p_t(\bfx)$ are concave for every $t\in[0,T]$, the sequence $(\eta^p)_{p\geq 0}$ is non-decreasing, and such that for any $p\geq 0$ and any $(t,\bfx)\in[0,T]\times\mathcal C^0$, the matrix $\sigma^{\eta^p}_t(\bfx)$ is an invertible matrix such that $(\sigma_t^{\eta^p})^{-1}(\bfx)f$ is uniformly bounded in $(t,\bfx)$, where
$$\sigma^{\eta^p}_t(\bfx):=\eta_t^p(\bfx) M_t+\sigma_t(\bfx).$$
$(v)$ The matrix $M_t\sigma^\top_t(\bfx)+\sigma_t(\bfx)M_t$ is symmetric negative, for every $(t,\bfx)\in[0,T]\times\mathcal C^0$.

\vspace{0.5em}
\noindent $(vi)$ The maps $U$, $\mu$ and $\sigma$ are such that $U$ is concave and for every $n\geq 1$, and $0\leq k\leq n-1,$ for every $(t,\bfx)\in[0,T]\times\mathcal C^0$,every $\{(\alpha_{i,j},\beta_{i,j},\gamma_{i,j})\in\R^d\times\R_+^*\times\R^d,\ 1\leq i\leq n-k,\ 0\leq j\leq n-k\}$, and every $(x_{0:k},\tilde\bfx,\tilde\bfy,\lambda)\in\left(\R^d\right)^{k+1}\times\mathcal C^t\times\mathcal C^t\times[0,1]$, we have
\begin{align*}
&U\Big(i_{n}\Big(x_{0:k-1},{\mathbf w}^\lambda(\tilde\bfx,\tilde\bfy;0,0),\sum_{i=0}^{1}{\mathbf w}^\lambda(\tilde\bfx,\tilde\bfy;i,1),\dots,\sum_{i=0}^{n-k}{\mathbf w}^\lambda(\tilde\bfx,\tilde\bfy;i,n-k)\Big)\Big)\\
&\geq U\Big(i_{n}\Big(x_{0:k-1},{\mathbf z}^\lambda(\tilde\bfx,\tilde\bfy;0,0),\sum_{i=0}^{1}{\mathbf z}^\lambda(\tilde\bfx,\tilde\bfy;i,1),\dots,\sum_{i=0}^{n-k}{\mathbf z}^\lambda(\tilde\bfx,\tilde\bfy;i,n-k)\Big)\Big)\end{align*}
where
\begin{align*}
{\mathbf w}^\lambda(\tilde\bfx,\tilde\bfy;i,\ell):=&\ \alpha_{i,\ell}+\beta_{i,\ell}\mu^{t,x}_{t^{t,n}_{k-1+i}}(\lambda\tilde\bfx+(1-\lambda)\tilde\bfy)+\Big(\eta_{t^{t,n}_{k-1+i}}^{t,\bfx}M_{t^{t,n}_{k-1+i}}+\sigma_{t^{t,n}_{k-1+i}}^{t,\bfx}\Big)(\lambda\tilde\bfx+(1-\lambda)\tilde\bfy)\gamma_{i,\ell},\\
{\mathbf z}^\lambda(\tilde\bfx,\tilde\bfy;i,\ell):=&\ \alpha_{i,\ell}+\beta_{i,\ell}\Big(\lambda\mu^{t,x}_{t^{t,n}_{k-1+i}}(\tilde\bfx)+(1-\lambda)\mu^{t,\bfx}_{t^{t,n}_{k-1+i}}(\tilde\bfy)\Big)+\gamma_{i,\ell}\lambda\Big(\eta_{t^{t,n}_{k-1+i}}^{t,\bfx}M_{t^{t,n}_{k-1+i}}+\sigma_{t^{t,n}_{k-1+i}}^{t,\bfx}\Big)(\tilde\bfx)\\
&+\gamma_{i,\ell}(1-\lambda)\Big(\eta_{t^{t,n}_{k-1+i}}^{t,\bfx}M_{t^{t,n}_{k-1+i}}+\sigma_{t^{t,n}_{k-1+i}}^{t,\bfx}\Big)(\tilde\bfy)
\end{align*}
\end{assumption} 

\begin{remark}\label{rem:assump}
This assumption deserves a certain number of comments. 
\begin{itemize}
\item $(iv')$ is here in order to ensure that the degenerate matrix $\sigma$ becomes invertible when it is suitably perturbed. Of course, our ultimate goal here is to assume that $\eta^p$ converges to $0$ and to approximate the solution of our problem with degenerate diffusion as the corresponding limit.
We also emphasize that this assumption implies in particular that for any $(t,\bfx)\in[0,T]\times\mathcal C^0$ and any $p\geq p'$
$$\sigma^{\eta^p}_t(\bfx)-\sigma^{\eta^{p'}}_t(\bfx)=(\eta_t^p(\bfx)-\eta^{p'}_t(\bfx))M_t,$$
which is a symmetric positive matrix. Hence, the sequence $\sigma^{\eta^p}$ is non-decreasing for the usual order on symmetric positive matrices.

\item $(v)$ and $(vi)$ are actually here mainly so that the results of Lemma \ref{lem.operator} below hold for a certain function $f$ involving $U$ (see the proof of Proposition \ref{prop:dege} below). They take a particularly complicated form for two reasons. First, our setting is fully non-Markovian, and second, it is also multidimensional. Indeed, as can be checked directly, if $d=1$ and $\bfx\longmapsto\sigma_t(\bfx)$ is linear, then we only need to assume that $\mu$ is concave and $f$ non-decreasing for \reff{eq:f} below to hold. Similarly, $(vi)$ is somehow a concavity assumption on $U$, $\mu$ and $\sigma$. Indeed, if again $d=1$ and if $U$ were Markovian, then a sufficient condition for $(vi)$ to hold is that $U$ is non-decreasing, $\mu$ is concave and $\sigma$ and $\eta$ are linear. 
\end{itemize}
\end{remark}

Our strategy of proof here is to start by obtaining a monotonicity result, with respect to the parameter $p$, for the solution of our control problem with diffusion coefficient $\sigma^p$. Such a result will be based on convex order type arguments. More precisely, we follow the strategy outlined by Pag\`es \cite{pages2014convex} and start by proving the result in a discrete--time setting, this is Proposition \ref{prop:dege}, which can then be extended to continuous-time through weak convergence arguments. Though the strategy of proof is the same as in \cite{pages2014convex}, our proofs are much more involved mainly due to the fact that, unlike in \cite{pages2014convex}, our framework is fully non--Markovian and multidimensional. We refer the reader to the appendix for details.
   
   \vspace{0.5em}
\noindent  Our main result is then that with degenerate volatility, the singular stochastic control problem can be represented as an infimum of solution of constrained BSDEs.
   \begin{theorem}\label{th:rep}
    Let Assumptions \ref{assump:U}, \ref{assump:exist}, \ref{assump:exist2} and \ref{ass: mu sigma lip2} hold, with $\sigma^p$ instead of $\sigma$, and assume in addition that 
    $$\sup_{(t,\bfx)\in[0,T]\times\mathcal C^0}\abs{\eta_t^p(\bfx)}\underset{p\rightarrow+\infty}{\longrightarrow}0.$$
    Then, we have
    $$v(t,\bfx)=\lim_{p\rightarrow +\infty}\uparrow v^p(t,\bfx)=\sup_{p>0}v^p(t,\bfx)=\sup_{p>0}\Yc^{t,\bfx,p}_t=\sup_{p>0}\Yw_t^{t,\bfx,p},$$
    where $\Yc^{t,\bfx,p}$ and $\Yw^{t,\bfx,p}$ are defined as $\Yc^{t,\bfx}$ and $\Yw^{t,\bfx}$ with $\sigma^p$ instead of $\sigma$.
   \end{theorem}
   
   \proof
   Since $\sigma^p$ satisfies all the required assumptions, by Theorem \ref{th:main} and Proposition \ref{prop:decreasing}, the only equality that we have to prove is the first one. But it is a simple consequence of classical estimates for SDEs and the uniform convergence we have assumed for $\eta^p$.
   \ep
   
   \begin{remark}
   The representation we have just obtained involves a supremum of solutions of constrained BSDEs. Formally speaking, such an object is close in spirit to so--called constrained second order BSDEs, as introduced by Fabre in her PhD thesis \cite{fabre2012some}. Indeed, the supremum over $p$ could be seen as a supremum over a family of probability measures, such that under these measures the canonical process has the same law as a continuous martingale whose quadratic variation has density $\sigma^p(\sigma^p)^\top$. To prove such a relationship rigorously is a very interesting problem, which however falls outside the scope of this paper.
   \end{remark}
   \subsection{Utility maximisation with transaction costs for non--Markovian dynamics}
        \subsubsection{The setup}
Fix $d=3$. We consider here a financial market, where the uncertainty will be driven by the third component of the underlying Brownian motion, namely for any $(t,\bfx)\in[0,T]\times \mathcal C^0$
$$X_s^{t,\bfx,3}:=\bfx^3_t+B_s^{t,3},\; s\in[t,T].$$
The market contains two traded assets. One is a non--risky asset, corresponding to a bank account, and whose value is given, for any $(t,\bfx)\in[0,T]\times \mathcal C^0$, by
$$\mathcal S_s^{t,\bfx}=\bfx^3_t+\int_t^sr_u^{t,\bfx^3}\big(X_u^{t,\bfx,3}\big)\mathcal S_u^{t,\bfx}du,\ s\in[t,T],$$
where $r$ is a bounded progressively measurable map, representing a non--Markovian interest rate on the market.

\vspace{0.5em}
The second asset is a risky one, and its value is given, for any $(t,\bfx)\in[0,T]\times \mathcal C^0$, by
                            $$S_s^{t,\bfx}=\int_t^sS_u^{t,\bfx}m^{t,\bfx^3}_u\big(X_u^{t,\bfx,3}\big)du+\int_t^sS_u^{t,\bfx}\Sigma^{t,\bfx^3}_u\big(X_u^{t,\bfx,3}\big)dB^{t,3}_u,\ s\in[t,T],$$
where $m$ and $\Sigma$ are bounded progressively measurable maps, representing respectively the drift and the volatility of the returns of the asset $S^{t,\bfx}$. We assume that these two maps are also Lipschitz--continuous in the sense of Assumption \ref{ass: mu sigma lip2}(iii).

\vspace{0.5em}
An investor has access to the market, but every transaction that he makes which involves the risky asset incurs a cost proportional to the number of asset bought or sold, given by a fixed constant $\lambda>0$ (see \cite{davis1990portfolio,shreve1994optimal} for more details and references). For any $t\in[0,T]$, we keep track of his transactions after time $t$ through the process $K\in \mathcal U_{\rm sing}^t$, in the sense that $K^1$ and $K^2$ record respectively the transactions from the risky to the non--risky asset and from the non--risky to the risky asset. Because of Proposition \ref{prop:simp}, we immediately restrict our attention to processes $K\in\overline{\mathcal U}^t$. Then, for any $\nu\in\mathcal U_t$, if $X^{t,\bfx,\nu,1}$ represents the total amount of money invested in the non--risky asset by the investor since time $t$, and $X^{\bfx,\nu,2}$ the total amount of money invested in the risky asset, we have immediately
$$\begin{cases}
               \displaystyle X_s^{t,\bfx,\nu,1}=\bfx_t^1+\int_t^sr^{t,\bfx^3}_u\left(X_s^{t,\bfx,\nu,3}\right)X_u^{t,\bfx,\nu,1}du+\int_t^s\left(\nu^1_u-(1+\lambda)\nu^2_u\right)du,\\[0.8em]
               \displaystyle X_s^{t,\bfx,\nu,2}=\bfx_t^2+\int_t^sX_u^{t,\bfx,\nu,2}\left(m^{t,\bfx ^3}_u\left(X_u^{t,\bfx,\nu,3}\right)du+\Sigma^{t,\bfx^3}_u\left(X_u^{t,\bfx,\nu,3}\right)dB^{t,3}_u\right)+\int_t^s\left(\nu^2_u-(1+\lambda)\nu^1_u\right)du.
                \end{cases}$$
Under this form, our general parameters are given by
 $$
            f:=\begin{pmatrix}
              1 & -(1+\lambda) &0\\
              -(1+\lambda) &1&0\\
              0&0&0
            \end{pmatrix},\;  \mu_t(\bfx) :=\begin{pmatrix} r_t(\bfx^3)\bfx^1_t\\ m_t(\bfx^3)\bfx^2_t\\ 0 \end{pmatrix},\ 
            \sigma_t(\bfx):=
                \begin{pmatrix} 
                    0& 0 &0\\
                    0 & 0 & \Sigma_t(\bfx^3)\bfx^2_t\\
                    0&0&1
                \end{pmatrix}.
        $$
        The goal of the investor is to maximise the utility of his terminal wealth, expressed in terms of liquidation value (meaning that assets held at time $T$ must be liquidated, which incurs a final transaction cost), and given by the following map $U$
$$U(\bfx)=:\mathcal U(\bfx^3,\ell(\bfx^1_T,\bfx^2_T)),$$
where the so--called liquidation function $\ell$ is defined by
$$\ell(x,y):=x+\frac{y^+}{1+\lambda}-(1+\lambda)y^-, \ (x,y)\in\R^2,$$
and the map $\mathcal U$ is assumed to be a (random) utility function, which is increasing and strictly concave with respect to its second variable, as well as locally Lipschitz--continuous with polynomial growth, so that Assumption \ref{assump:U} is satisfied.

\vspace{0.5em}
The value function of the investor is thus given by
\begin{equation}\label{eq:transaction}v(t,\bfx)=\underset{\nu\in\mathcal U^t}{\sup}\mathbb E^{\P_0^t}\big[\mathcal U\big(\bfx^3\otimes_tX^{t,\bfx,\nu,3},\ell\big(X^{t,\bfx,\nu,1}_T,X^{t,\bfx,\nu,2}_T\big)\big)\big].
\end{equation}
Since $\sigma$ is clearly not invertible, we will perturb the above dynamics by adding some noise in the two above equations. More precisely, we define for any $p>0$ 
$$\begin{cases}
               \displaystyle X_s^{t,\bfx,\nu,p,1}=\bfx_t^1+\int_t^sr^{t,\bfx^3}_u\left(X_s^{t,\bfx,\nu,p,3}\right)X_u^{t,\bfx,\nu,p,1}du-\frac1pB^{t,1}_s+\int_t^s\left(\nu^1_u-(1+\lambda)\nu^2_u\right)du,\\[0.8em]
               \displaystyle X_s^{t,\bfx,\nu,p,2}=\bfx_t^2+\int_t^sX_u^{t,\bfx,\nu,p,2}\left(m^{t,\bfx ^3}_u\left(X_u^{t,\bfx,\nu,p,3}\right)du+\Sigma^{t,\bfx^3}_u\left(X_u^{t,\bfx,\nu,p,3}\right)dB^{t,3}_u\right)+\int_t^s\left(\nu^2_u-(1+\lambda)\nu^1_u\right)du\\[0.5em]
               \hspace{5em}               \displaystyle-\frac1pB^{t,2}_s,\\[0.8em]
                 \displaystyle      X_s^{t,\bfx,\nu,p,3}= X_s^{t,\bfx,3}.
                \end{cases}$$
Still with our notations, this corresponds to
        $$
          \eta_t^p(\bfx)=-\frac1p,\  M_t:=\begin{pmatrix} 
                    1& 0 &0\\
                    0 & 1 & 0\\
                    0&0&0
                \end{pmatrix}.
        $$
     Notice then that
        $$
                (\sigma^p_t)^{-1}(\bfx)=
                \begin{pmatrix}
                    -p & 0 &0\\
                    0 & -p & p\sigma^S_t(\bfx^3)\bfx^2_t\\
                    0&0&1
                \end{pmatrix},\ M_t\sigma_t(\bfx)=\begin{pmatrix} 
                    0& 0 &0\\
                    0 & 0 & 0\\
                    0&0&0
                \end{pmatrix}, \text{ and }(\sigma^p_t)^{-1}(\bfx)f= - pf,$$
                so that Assumption \ref{ass: mu sigma lip2}(iv'),(v) is satisfied.

\subsubsection{The result}
\noindent We will actually use the result proved in Theorem \ref{th:rep} conditionally on $X^{t,\bfx,\nu,p,3}$ (that is to say that we consider conditional expectations with respect to $\sigma(X^{t,\bfx,\nu,p,3}_s,\ t\leq s\leq T)$ instead of simple expectations). It can be checked readily that all our arguments still go through in this case. Moreover, the drifts and volatility in the dynamics of $X^{t,\bfx,\nu,p,1}$ and $X^{t,\bfx,\nu,p,2}$ then become (conditionally) linear. Then, remembering Remark \ref{rem:assump} above, we know that (conditionally), Assumption \ref{ass: mu sigma lip2}(vi) is also satisfied. Notice also that in this case the constraint can be read
$$-pf\Zc^{t,\bfx,p}_s\in \mathfrak K_s,\; dt\times d\P-a.e.,$$
which, by definition of $\mathfrak K$ and since $p>0$, is actually equivalent to $-f\Zc_s^{t,\bfx,p}\in \mathfrak K_s$. Our main result of thus section is 
 \begin{theorem}
The value function of the utility maximisation problem in finite horizon with proportional transaction costs \eqref{eq:transaction} satisfies
$$v(t,\bfx)= \underset{p>0}{\sup}\; \Yc^{t,\bfx,p}_t,$$
where $(\Yc^{t,\bfx,p},\Zc^{t,\bfx,p})$ is the minimal solution of the $Z-$constrained {\rm BSDE}
   \begin{align*}
            \Yc^{t,\bfx,p}_\cdot\ge \mathcal U\big(\bfx^3\otimes_t\bfX^{t,\bfx,p,3},\ell\big(\bfX^{t,\bfx,p,1}_T,\bfX^{t,\bfx,p,2}_T\big)\big)-\int_\cdot^T \Zc^{t,\bfx,p}_s\cdot dB^t_s,\ \P^t_0-a.s.,\;
            -f\Zc^{t,\bfx,p}_s\in \mathfrak K_s, \ ds\otimes d\P^t_0-a.e.
        \end{align*}

 \end{theorem}

\begin{remark}
Notice that the BSDE representation of the value function depends on the parameter $p$ only through the terminal condition. Standard estimates imply that this terminal condition converges in $L^p(\P_0^t)$ to $U\big(\bfx^3\otimes_t\bfX^{t,\bfx,3},\ell\big(\bfX^{t,\bfx,1}_T,\bfX^{t,\bfx,2}_T\big)\big)$. Using then stability results for constrained BSDEs, similar to the ones proved in the next section, we reasonably expect that $\Yc^{t,\bfx,p}_t$ will converge to $Y^{t,\bfx}$, where $(Y^{t,\bfx},Z^{t,\bfx})$ is the minimal solution of the $Z-$constrained {\rm BSDE}
   \begin{align*}
            Y^{t,\bfx}_\cdot\ge \mathcal U\big(\bfx^3\otimes_t\bfX^{t,\bfx,3},\ell\big(\bfX^{t,\bfx,1}_T,\bfX^{t,\bfx,2}_T\big)\big)-\int_\cdot^T Z^{t,\bfx}_s\cdot dB^t_s,\ \P^t_0-a.s.,\;
            -fZ^{t,\bfx}_s\in \mathfrak K_s, \ ds\otimes d\P^t_0-a.e.
        \end{align*}
Proving this rigorously falls however outside the scope of this paper.
\end{remark}

\noindent As far as we know, such a result is completely new in the literature, even in the Markovian case. Moreover, as pointed out in the recent paper \cite{kallsen2013portfolio}, the non--Markovian case has actually never been studied using stochastic control and PDE tools, the only approach in the literature being convex duality. We thus believe that our approach achieves a first step allowing to tackle this difficult problem. Let us nonetheless point out a gap in our approach. If we wanted to cover completely the problem of transaction costs, we should have added state constraints in our original stochastic control problem. Indeed, those are inherent to the problem of transaction costs, in order to avoid bankruptcy issues (though this is actually a lesser issue when the time horizon is finite, as in our case). We have chosen not to do so so as not to complicate even more our arguments, but we believe that they could be also used in this setting, albeit with possibly important modifications. In particular, the full dynamic programming principle that we used does not seem to be proved in such a general framework in the literature, when state constraints are present (see however \cite{bouchard2012weak,bouchard2015stochastic}), and it is not completely clear in which sense the equality between $v^{\rm sing}$ and $v$ will then hold.

\section{Applications: DPP and regularity for singular stochastic control}\label{Sec6}
\subsection{Dynamic programming principle}
Notice that in all our study, we never actually proved that the dynamic programming principle actually held for the singular stochastic control problem defining $v$ (or $v^{\rm sing}$). However, it is an easy consequence of our main result.

\begin{theorem}
For any $(t,\bfx)\in[0,T]\times\mathcal C^0$, for any $\tau\in\Tc^t$ and any $\theta\in\Tc^{t,\bfx}$, we have
                \begin{align*}
                     v(t,\bfx)&=\underset{\nu\in\Uc^{t}}{\sup}\E^{\P_0^t}\left[v(\tau,\bfx\otimes_tX^{t,\bfx,\nu})\right]=\underset{\nu\in\Uc^{t}}{\sup}\E^{\P_\nu^{t,\bfx}}\left[v(\theta,\bfx\otimes_tB^{t,\bfx_t})\right].
                                            \end{align*}
\end{theorem}

\proof
It is an immediate consequence of the dynamic programming principle satisfied by the penalised BSDEs \eqref{eq: ppde for penalized bsde dpp} and the convergence of penalised BSDEs to the minimal solution of the constrained BSDE.
\ep
\subsection{Regularity results}
In this section, we show how our representation can help to obtain {\it a priori} regularity results for the value function of singular stochastic control problems. Such results in that level of generality are, as far as we know, the first available in the literature. 

\vspace{0.5em}
The main idea of the proof is that as soon as one knows that the value function of the singular control problem is associated to a constrained BSDE, one can use the fact that such BSDEs are actually linked to another different singular stochastic control problem, which is actually simpler to study. Such a representation is not new and was already the crux of the arguments of Cvitani\'c, Karatzas and Soner \cite{cvitanic1998backward}. It has also been used very recently in \cite{bouchard2014regularity} to obtain the first regularity results in the literature for constrained BSDEs. For the sake of simplicity, and since this is not the heart of our article, we will concentrate on the Markovian set--up for this application, and leave the more general case to future research\footnote{We would like to point out the reader to the recent work in preparation \cite{bouchard2017prep} which will actually extend the results of \cite{bouchard2014regularity} to the non-Markovian case}.

\vspace{0.5em}
Let us define the map $\delta:[0,T]\times\mathbb R^d\longrightarrow \R_+$ such that for any $t\in[0,T]$, $\delta_t(\cdot)$ is the so-called support function of the set $K_t$, that is to say
$$\delta_t(u):=\sup\left\{k\cdot u,\ k\in K_t\right\}.$$
Notice that since the zero vector in $\mathbb R^d$ belongs to $K_t$ for any $t\in[0,T]$, it is clear that $\delta$ is non--negative. 
This section requires requires the following additional assumptions.
\begin{assumption}\label{assump:delta}
$(i)$ The maps $U$, $\mu$ and $\sigma$ are Markovian, that is to say, abusing notations slightly
$$U(\bfx)=U(\bfx_T),\ \mu_t(\bfx)=\mu_t(\bfx_t),\ \sigma_t(\bfx)=\sigma_t(\bfx_t),\, \text{for any $\bfx\in\mathcal C^0$}.$$
$(ii)$ the map $t\longmapsto f_t$ does not depend on $t$, and thus $\delta$ as well.

\vspace{0.5em}
$(iii)$ If one defines the face--lift of $U$ by
$$\widehat U(x):=\sup_{u\in\mathbb R^d}\left\{U(x+f_Tu)-\delta(u)\right\},\ x\in\mathbb R^d,$$
then we have for some constant $C>0$ and any $(x,x')\in\mathbb R^d\times\mathbb R^d$
$$\abs{\widehat U(x)-\widehat U(x')}\leq C\No{x-x'}.$$
\end{assumption}

The main result of this section is
\begin{theorem}\label{th:regularity}
Let Assumptions \ref{ass: mu sigma lip}, \ref{assump:U}, \ref{assump:exist}, \ref{assump:exist2} and \ref{assump:delta} hold. Then, there is a constant $C>0$ such that for any $(t,t',\bfx,\bfx')\in[0,T]\times[t,T]\times\mathcal C^0\times\mathcal C^0$
$$\abs{v(t,\bfx)-v(t,\bfx')}\leq C\No{\bfx_t-\bfx'_t},\  \abs{v(t,\bfx)-\E^{\P_0^t}[v(t',\bfx)]}\leq C\big(1+\No{\bfx_t}\big)(t'-t)^{1/2}.$$
\end{theorem}

The remaining of this section is dedicated to the proof of this result. We shall make a strong use of the connection with constrained BSDEs established before. 
\subsubsection{Another singular control problem}\label{sec:delta}
For any $t\in[0,T]$, let us consider the following set of controls
\begin{align*}
\mathcal V^t_{\rm b}:=&\Big\{(u_s)_{t\leq s\leq T},\ \text{which are $\R^d-$valued, $\F^t-$predictable and bounded.}\Big\}.
\end{align*}

For any $(t,x)\in[0,T]\x\mathbb R^d$, we define
$$Y^{x}_t:=\sup_{u\in\mathcal V^t_{\rm b}}\mathbb E^{\mathbb P_0^t}\left[U(X^{t,x,u}_T)-\int_t^T\delta(u_s)ds\right],$$
where $X^{t,x,u}$ is the unique strong solution on $(\Omega^t,\mathcal F_T^{t,o},\mathbb P_0^t)$ of the following SDE
    \begin{align*}
        \displaystyle
            X^{t,x,u}&=x+\int_t^{\cdot} \mu_s\big(X_s^{t,x,u}\big)ds+\int_t^{\cdot}f_su_sds +\int_t^{\cdot} \sigma_s\big(X_s^{t,x,u}\big)dB^t_s\;,
                \;\P_0^t\mbox{-a.s.}\nonumber
    \end{align*}

This value function is always well-defined since $u$ is bounded and $\delta$ is non--negative. Our first step is to show that one can actually replace the map $U$ above by its facelift. It is a version of Proposition 3.1 of \cite{bouchard2014regularity} for our setting. The proof is postponed to the appendix.
\begin{lemma}\label{lemma:face}
Let Assumptions \ref{ass: mu sigma lip}, \ref{assump:U} and \ref{assump:delta} hold. Then, for any $t<T$, we have
$$Y^{x}_t=\sup_{u\in\mathcal V^t_{\rm b}}\mathbb E^{\mathbb P_0^t}\left[\widehat U(X_T^{t,x,u})-\int_t^T\delta(u_s)ds\right].$$
\end{lemma}

\vspace{0.5em}
The next result is Proposition 3.3 of \cite{bouchard2014regularity} in our framework
\begin{lemma}\label{lemma:facelift}
Let Assumptions \ref{ass: mu sigma lip}, \ref{assump:U} and \ref{assump:delta} hold. We have for any $(t,x)\in[0,T)\times\mathbb R^d$
$$Y_t^{x}=\underset{\iota\in \mathbb L_\infty(\Fc_t)}{\rm essup}\left\{Y_t^{x+f\iota}- \delta(\iota)\right\},\ a.s.,$$
where $\mathbb L_\infty(\Fc_t)$ is the set of $\R^d$-valued, bounded and $\Fc_t$-measurable random variables.
\end{lemma}

\vspace{0.5em}
We can now give our main result of this section.

\begin{proposition}\label{prop:prop}
Let Assumptions \ref{ass: mu sigma lip}, \ref{assump:U} and \ref{assump:delta} hold. Then, there is some constant $C>0$ such that for any $0\leq t\leq s< T$, any $(x,x')\in\mathbb R^d\times\mathbb R^d$
$$\abs{Y_t^x-Y_t^{x'}}\leq C\No{x-x'},\ \abs{Y_t^{x}-\mathbb E^{\mathbb P_0^t}[Y_{s}^x]}\leq C(1+\No{x})(s-t)^{\frac 12} .$$
\end{proposition}

\subsubsection{Weak formulation and the main result}
 For any $(t,x)\in[0,T]\x\mathbb R^d$ and $u\in\Vc^t_{\rm b}$, we now define the following $\P_0^t-$equivalent measure
    $$
        \frac{d\P^{t,x,u}}{d\P_0^t}=\Ec\left( \int_t^\cdot (\sigma_s)^{-1}\left( X^{t,x,0} \right)fu_s\cdot dB^t_s \right).
    $$
    The weak formulation of the control problem is defined as
    \begin{equation*}
    Y^{{\rm w},x}_t:=\underset{u\in\mathcal V^t_{\rm b}}{\sup} \mathbb E^{\mathbb P^{t,x,u}}\left[U(X^{t,x,0})-\int_t^T\delta(u_s)\right],\ \text{for any $(t,x)\in[0,T]\times\mathbb R^d$}.
    \end{equation*}
The following proposition is a simple consequence of Remark 3.8 and Theorem 4.5 of \cite{karoui2013capacities2}.
\begin{proposition}\label{prop:w}
For any $(t,x)\in[0,T]\times\mathbb R^d$, we have $Y^x_t= Y^{{\rm w},x}_t$.
\end{proposition}
We can now proceed to the 
\proof[Proof of Theorem \ref{th:regularity}]
By Theorem 4.1 of \cite{bouchard2014regularity}, we have for any $(t,\bfx)\in[0,T]\times\mathcal C^0$, $Y^{{\rm w},\bfx(t)}_t=\mathcal Y^{t,\bfx}_t$ defined as the first component of the constrained BSDE \eqref{eq:bsde1}--\eqref{eq:bsde2}. By Propositions \ref{prop:prop} and \ref{prop:w} we then deduce that there is a constant $C>0$ such that for any $(t,t',\bfx,\bfx')\in[0,T]\times[t,T]\times\mathcal C^0\times\mathcal C^0$
$$\abs{\mathcal Y^{t,\bfx}_t-\mathcal Y^{t,\bfx'}_t}\leq C\No{x_t-x'_t},\ \abs{\mathcal Y^{t,\bfx}_t-\mathbb E^{\mathbb P_0^t}[\mathcal Y^{t',\bfx}_{t'}]}\leq C(1+\No{\bfx})(t'-t)^{\frac 12} .$$
It then suffices to apply Theorem \ref{th:main}.
\ep
  \bibliographystyle{plain}
 \small
\bibliography{bibliographyDylan}
\appendix
\section{Technical proofs for Section \ref{Sec2}}
\proof[Proof of Proposition \ref{prop:simp}]
$(i)$ First of all, it is a classical result that $\mathcal U^t$ is dense in $\mathcal U_{\rm sing}^t$ in the sense that for any $K\in\mathcal U_{\rm sing}^t$, there is some sequence $(\nu^n)_{n\geq 0}\subset\mathcal U^t$ such that
 \begin{equation}\label{eq:conv}\mathbb E^{\P_0^t}\Bigg[\underset{t\leq s\leq T}{\sup}\bigg\|K_s-\int_t^s\nu^n_rdr\bigg\|^2\Bigg]\underset{n\rightarrow+\infty}{\longrightarrow}0.\end{equation}

It is also clear from Assumption \ref{assump:U}, \eqref{eq: estimate SDE growth} and \eqref{eq: estimate SDE Delta} with $\bfx=\bfx'$ that for any $(t,\bfx;K,K')\in[0,T]\x\mathcal C^0\x(\Uc^t_{\rm sing})^2$, we have for some constant $C_0>0$ which may vary from line to line
\begin{align*}
&\abs{\mathbb E^{\P_0^t}\big[U\big(\bfx\otimes_tX^{t,\bfx,K}\big)\big]-\mathbb E^{\P_0^t}\big[U\big(\bfx\otimes_tX^{t,\bfx,K'}\big)\big]}\\
&\leq C_0\mathbb E^{\P_0^t}\Big[\big\|\bfx\otimes_t\big(X^{t,\bfx,K}-X^{t,\bfx,K'}\big)\big\|_{\infty,T}^2\Big]^{\frac12}\Big(1+\mathbb E^{\mathbb P_0^t}\Big[\big\|\bfx\otimes_tX^{t,\bfx,K}\big\|_{\infty,T}^{2r}\Big]^{\frac12}+\mathbb E^{\mathbb P_0^t}\Big[\big\|\bfx\otimes_tX^{t,\bfx,K'}\big\|_{\infty,T}^{2r}\Big]^{\frac12}\Big)\\
&\leq C_0\E^{\P_0^t}\Bigg[\bigg\|\int_t^Td\big(K_s-K'_s\big)\bigg\|^{2}\Bigg]^{\frac12}\Bigg(1+\No{\bfx}_{\infty,t}^r +\E^{\P_0^t}\Bigg[\bigg\|\int_t^TdK_s\bigg\|^{2r}\Bigg]^{\frac12}+\E^{\P_0^t}\Bigg[\bigg\|\int_t^TdK'_s\bigg\|^{2r}\Bigg]^{\frac12}\Bigg).
\end{align*}

We thus deduce immediately that the map $K\longmapsto \mathbb E^{\P_0^t}\left[U\left(\bfx\otimes_tX^{t,\bfx,K}\right)\right]$ is continuous with respect to the convergence in \eqref{eq:conv}. Hence the first result.

\vspace{0.5em}
$(ii)$ It is a classical result and we refer the reader to Proposition 4.1 in \cite{bouchard2014regularity} or Theorem 4.5 of \cite{karoui2013capacities2}.

\vspace{0.5em}
$(iii)$ Notice that by definition, $\overline{\F^{\bfX^{t,\bfx}}}^{\mathbb P_0^t}$ satisfies the Blumenthal $0-1$ law as well as the predictable martingale representation property. This implies that the process below is a Brownian motion on $(\bfL{t,\bfx_t},\Fc^{t,\bfx_t,o}_T,\F^{t,\bfx_t},\P^{t,\bfx}_{0} )$
    $$W^{t,\bfx}:=\int_t^\cdot\big(\sigma^{t,\bfx}_s\big)^{-1}\big(B^{t,\bfx_t}\big)\big(dB^{t,\bfx_t}_s-\mu^{t,\bfx}_s\big(B^{t,\bfx_t}\big)ds\big),\; \mathbb P^{t,\bfx}_0\mbox{-a.s.}$$
    Indeed, by definition of $\P^{t,\bfx}_0$, we know that
    \begin{equation}\label{eq:lawlaw}
    \text{the law of $B^{t,\bfx_t}$ under $\mathbb P_0^{t,\bfx}$ $=$ the law of $\bfX^{t,\bfx}$ under $\P_0^t$}.\end{equation}
    Hence, since we do have
    $$B^{t}=\int_t^{\cdot}\big(\sigma^{t,\bfx}_s\big)^{-1}\big(\bfX^{t,\bfx}\big)\big(d\bfX^{t,\bfx}_s-\mu^{t,\bfx}_s\big(\bfX^{t,\bfx}\big)ds\big),\: \mathbb P^{t}_0{\mbox{-a.s.}},$$
    the result follows immediately by the L\'evy characterisation of Brownian motion and \eqref{eq:filt}. Moreover, we clearly also have that
    \begin{equation}\label{eq:W}
    \text{the law of $W^{t,\bfx}$ under $\mathbb P_0^{t,\bfx}$ $=$ the law of $B^t$ under $\P_0^t$}.
    \end{equation}

    Next, for any $\nu\in\Uc^t$, we define the following probability measure $\P^{t,\bfx}_\nu$ on $(\Lambda^{t,\bfx_t},\Fc^{t,\bfx_t,o}_T)$
    $$
        \frac{d\P^{t,\bfx}_\nu}{d\P_0^{t,\bfx}}:=\Ec\left( \int_t^\cdot \big(\sigma^{t,\bfx}_s\big)^{-1}\big( B^{t,\bfx_t} \big)f_s\nu_s\big(W^{t,\bfx}\big)dW^{t,\bfx}_s \right)_T,
    $$
    where it is understood that we interpret $\nu$ as a (Borel) map from $\mathcal C^t$ to $\mathbb R^d$. We claim that for all $t\in[0,T]$
\begin{equation*}
\E^{\P^{t,\bfx,\nu}}\big[ U^{t,\bfx}\big( \bfX^{t,\bfx} \big) \big]=\E^{\P^{t,\bfx}_\nu}\big[ U^{t,\bfx}\big( B^{t,\bfx} \big) \big],\; \text{for every $\nu\in\Uc^t$.}
\end{equation*}
Indeed, we have, using successively the definition of $\P^{t,\bfx}_\nu$, \eqref{eq:lawlaw}, \eqref{eq:W} and the definition of $\P^{t,\bfx,\nu}$
\begin{align*}
\E^{\P^{t,\bfx}_\nu}\big[U^{t,\bfx}\big(B^{t,\bfx_t}\big)\big]&=\E^{\P_0^{t,\bfx}}\bigg[\Ec\bigg( \int_t^\cdot \big(\sigma^{t,\bfx}_s\big)^{-1}\big( B^{t,\bfx_t} \big)f_s\nu_s\big(W^{t,\bfx}\big)dW^{t,\bfx}_s \bigg)_TU^{t,\bfx}\big(B^{t,\bfx_t}\big)\bigg]\\
&=\E^{\P_0^{t}}\bigg[\Ec\bigg( \int_t^\cdot \big(\sigma^{t,\bfx}_s\big)^{-1}\big( \bfX^{t,\bfx} \big)f_s\nu_s\big(B^{t}\big)dB^{t}_s \bigg)_TU^{t,\bfx}\big(\bfX^{t,\bfx}\big)\bigg]\\
&=\E^{\P^{t,\bfx,\nu}}\big[U^{t,\bfx}\big(\bfX^{t,\bfx}\big)\big].
\end{align*}
The result is then immediate.
\ep

\vspace{0.5em}
 \proof[Proof of Lemma \ref{lemma:conv}]
 It is clear that the sequence is non-decreasing, as the sequence of sets $\mathcal U^{.,n}$ is. Moreover, since the elements of $\mathcal U^t$ have by definition moments of any order, it is clear that $\cup_{n\geq 1}\mathcal U^{t,n}$ is dense in $\mathcal U^t$, in the sense that for any $\nu\in\mathcal U^t$, there exists a sequence $(\nu^m)_{m\geq 1}$ such that for any $m\geq 1$, $\nu^m\in \mathcal U^{t,m}$ and 
\begin{equation}\label{eq:conv2}\mathbb E^{\P_0^t}\left[\int_t^T\No{\nu_r-\nu^m_r}^2dr\right]\underset{m\rightarrow+\infty}{\longrightarrow}0.\end{equation}

By the same arguments as in the proof of Proposition \ref{prop:simp}, we deduce that
$$v(t,\bfx)=\underset{\nu\in\cup_{n\geq 1}\mathcal U^{t,n}}{\sup}\mathbb E^{\mathbb P^t_0}\left[U^{t,\bfx}\left(X^{t,\bfx,\nu}\right)\right]=\underset{n\rightarrow +\infty}{\lim} v^n(t,\bfx),$$
since the sets $\mathcal U^{t,n}$ are non-decreasing with respect to $n$.
 \ep
 
 \section{Technical proofs for Section \ref{Sec4}}
  \proof[Proof of Lemma \ref{lemma:utile}]
    
        Fix $(t,\bfx)\in[0,T]\x\mathcal C^0$ and $\tau\in\Tc^{t,\bfx}_+$ and denote by $(U,V)\in\mathfrak S^2_{t,\bfx}\times\mathfrak H^2_{t,\bfx}$ the unique solution of the following backward stochastic differential equation on $[t,\tau]$
        $$
            U_s=\tau-t-\int_s^\tau\bar n\abs{V_r}dr-\int_s^\tau V_r\cdot dW^{t,\bfx}_r,\ s\in[t,\tau],\ \mathbb P^{t,\bfx}_0-a.s.,
        $$
       where we remind the reader that in our context $\tau$ is $\Fc^{t,\bfx,o}_\tau$-measurable, and that we can always consider a $\mathbb P_0^{t,\bfx}$-version of $U$ (resp. $V$), which we still denote by $U$ (resp. $V$) for simplicity, and which is $\mathbb F^{t,\bfx_t,o}-$progressively measurable (resp. predictable).
       
       \vspace{0.5em}
        Let $\mu$ be an arbitrary $\F^{t,\bfx_t,o}-$predictable process satisfying $\abs\mu\le\bar n$, so that for all $z\in\R^d$ we have 
        $-\bar n\abs z\le \mu\cdot z$. 
        This implies in particular that $\Q^\mu\in\Mc^{t,\bfx,\bar n}$ with 
        $$
            d\Q^\mu:=\Ec\left(\int_t^T\mu_s\cdot dW_s^{t,\bfx}\right)d\P_0^{t,\bfx}.
        $$
        Hence, by standard a comparison result for BSDEs (see for instance Theorem 2.2 in \cite{el1997backward}) we have
        \begin{align*}
        U_t\le \E^{\Q^\mu}[\tau-t].
        \end{align*}
        Hence, the arbitrariness of $\mu$ implies that
        \beq\label{eq: inf tau positive ineq 1}
            U_t\le \underline\Ec^n_t[\tau-t].
        \eeq
        On the other hand, let $\nu^n$ be defined such that $\nu^n\cdot V=-\bar n\abs V$ and observe that $\abs{\nu^n}\le\bar n$.
        Then we have $\Q^n\in\Mc^{t,\bfx,\bar n}$ with
        $$
            d\Q^n:=\Ec\left(\int_t^T\nu^n_s\cdot dW_s^{t,\bfx}\right)d\P_0^{t,\bfx}.
        $$
        We hence have, using the fact that $U_t$ is a constant by the Blumenthal $0-1$ law
        $$
            U_t=\tau-t-\int_t^\tau V_s\cdot\left( dW^{t,\bfx}_s - \nu^n_s ds \right)=\E^{\Q^n}[\tau-t]\ge \underline\Ec^n_t[\tau-t],
        $$
        by definition of $\underline\Ec^n_t$.
        The last inequality together with \reff{eq: inf tau positive ineq 1} gives that
        $$
            U_t=\underline\Ec^n_t[\tau-t].
        $$
        Since $\tau>t$, $\P^{t,\bfx}_0$-a.s., the result follows by strict comparison (see again Theorem 2.2 in \cite{el1997backward}).
    \ep

\vspace{0.5em}
  \proof[Proof of Proposition \ref{prop: ppde for app vf}]
        $(i)$ The case of $v^n$.

\vspace{0.5em}
We proceed along the lines of \cite[Proof of Proposition 4.4]{ekren2016viscosity}, splitting the proof into two steps.

\vspace{0.5em}
Fix $(t,\bfx)\in[0,T]\x\mathcal C^0$ for the remainder of this part.

\vspace{0.5em}

            \step1{We show that $v^{n}\in C^0([0,T]\times\bfL{}^{0,\bfx_0})$ and satisfies the dynamic programming principle,
                for any $\tau\in\Tc^t$ and any $\theta\in\Tc^{t,\bfx}$
                \begin{equation}\label{eq: ppde for app vf dpp}
                     v^n(t,\bfx)=\underset{\nu\in\Uc^{t,n}}{\sup}\E^{\P_0^t}\left[v^n(\tau,\bfx\otimes_tX^{t,\bfx,\nu})\right]=\underset{\nu\in\Uc^{t,n}}{\sup}\E^{\P_\nu^{t,\bfx}}\left[v^n(\theta,\bfx\otimes_tB^{t,\bfx_t})\right]
\end{equation}}

The dynamic programming result is actually classical since the controls $\nu$ that we consider here take values in a compact subset of $\mathbb R^d$. We refer the reader to Proposition 2.5 and Theorem 3.3 in \cite{karoui2013capacities2}, where we emphasize that their proof of Theorem 3.3 can immediately be extended to the case of $\mu$ and $\sigma$ Lipschitz continuous with linear growth (instead of Lipschitz continuous and bounded), as they only require to have existence of a strong solution to the SDEs considered.

\vspace{0.5em}
We next show the continuity of $v^n$. For any $(t,t';\bfx,\bfx')\in[0,T]\times[t,T]\times\mathcal C^0\times\mathcal C^0$, we have
\begin{align*}
\abs{v^n(t,\bfx)-v^n(t',\bfx')}\leq \abs{v^n(t,\bfx)-v^n(t,\bfx')}+\abs{v^n(t,\bfx')-v^n(t',\bfx')}.
\end{align*}
We now estimate separately the two terms in the right-hand side above. We have first using Assumption \ref{assump:U}, \eqref{eq: estimate SDE growth}, \eqref{eq: estimate SDE Delta} and the fact that the controls $\nu\in\mathcal U^{t,n}$ are bounded by $n\sqrt{d}$
\begin{align*}
&\abs{v^n(t,\bfx)-v^n(t,\bfx')}\leq \underset{\nu\in\mathcal U^{t,n}}{\sup} \mathbb E^{\mathbb P_0^t}\Big[\big|U^{t,\bfx}\big(X^{t,\bfx,\nu}\big)-U^{t,\bfx'}\big(X^{t,\bfx',\nu}\big)\big|\Big]\\
&\leq  C\underset{\nu\in\mathcal U^{t,n}}{\sup} \mathbb E^{\mathbb P_0^t}\Big[\big\|\bfx\otimes_t\big(X^{t,\bfx,\nu}-X^{t,\bfx',\nu}\big)\big\|_{\infty,T}^2\Big]^{\frac12}\Big(1+\mathbb E^{\mathbb P_0^t}\Big[\big\|\bfx\otimes_tX^{t,\bfx,\nu}\big\|_{\infty,T}^{2r}\Big]^{\frac12}+\mathbb E^{\mathbb P_0^t}\Big[\big\|\bfx'\otimes_tX^{t,\bfx',\nu}\big\|_{\infty,T}^{2r}\Big]^{\frac12}\Big)\\
&\leq C_dn^{1\vee r}\No{\bfx-\bfx'}_{\infty,t}\left(1+\No{\bfx}_{\infty,t}^r+\No{\bfx'}_{\infty,t}^r\right),
\end{align*}
where the constant $C_d$ does not depend on $n$.

\vspace{0.5em}
Next, using \eqref{eq: ppde for app vf dpp} for $\tau=t'$, we compute, using the previous calculation and \eqref{eq: estimate SDE Delta initial}
\begin{align*}
&\abs{v^n(t,\bfx')-v^n(t',\bfx')}\leq \underset{\nu\in\mathcal U^{t,n}}{\sup}\E^{\mathbb P_0^t}\Big[\big|v^n(t',\bfx'\otimes_tX^{t,\bfx',\nu})-v^n(t',\bfx')\big|\Big]\\
&\leq C_dn^{1\vee r}\underset{\nu\in\mathcal U^{t,n}}{\sup}\E^{\mathbb P_0^t}\Big[\big\|\bfx'\otimes_tX^{t,\bfx',\nu}-\bfx'\big\|_{\infty,t'}^2\Big]^{\frac12}
\Big(1+\No{\bfx'}_{\infty,t'}^r+\E^{\mathbb P_0^t}\Big[\big\|\bfx'\otimes_tX^{t,\bfx',\nu}\big\|^{2r}_{\infty,t'}\Big]^{\frac12}\Big)\\
&\leq C_dn^{1\vee r +1+r}(t'-t)^{\frac12}\left(1+\No{\bfx'}^{r+1}_{\infty,t}+\No{\bfx'}^{r+1}_{\infty,t'}\right).
\end{align*}

By definition of $d_\infty$, we have thus obtained that
$$\abs{v^n(t,\bfx)-v^n(t',\bfx')}\leq C_dn^{1+r+1\vee r}d_\infty\left((t,\bfx),(t',\bfx')\right)\left(1+\No{\bfx}_{\infty,t}^{r+1}+\No{\bfx'}_{\infty,t'}^{r+1}\right),$$
which proves the continuity of $v^n$ with respect to $d_\infty$.

\vspace{0.5em}
            \step2{We show that $v^n$ is a viscosity sub--solution to {\rm PPDE} \reff{eq: penalized PPDE}.}

\vspace{0.5em}
                Assume to the contrary that there is $(t,\bfx;\vp)\in[0,T]\x\mathcal C^0\x\underline\Ac v^n(t,\bfx)$ s.t. for some $c>0$
                $$
                    -\Lc^{t,\bfx}\vp(t,\bfx^t)-n\rho(f_t^\top D\vp(t,\bfx^t))\ge2c>0.
                $$
                Without loss of generality, we may reduce $H$ in the definition of $\vp\in\underline\Ac v^n(t,\bfx)$
                so that by continuity of all the above maps, we obtain
                $$
                    -\Lc^{t,\bfx}\vp(s,B^{t,\bfx_t})-n\rho(f_s^\top D\vp(s,B^{t,\bfx_t}))\ge c,\ \mbox{on}\  [t,H],\ \mathbb P_0^{t,\bfx}-a.s.
                $$
                Furthermore, observe that for each $s\in[t,H]$
                $$
                    n\rho(f_s^\top D\vp(s,B^{t,\bfx_t}))=\sup_{u\in[0,n]^d}u\cdot (f_s^\top D\vp(s,B^{t,\bfx_t})),
                $$
                so that by definition of $\Uc^{t,n}$ we have for all $\nu\in\Uc^{t,n}$
                \beq\label{eq: control subsol}
                    -\Lc^{t,\bfx}\vp(s,B^{t,\bfx_t})-\nu_s\cdot f_s^\top D\vp(s,B^{t,\bfx_t})\ge c, \ \mbox{on}\ [t,H],\ \mathbb P_0^{t,\bfx}-a.s.
                \eeq
                Since $\vp$ is smooth on $[t,H]$, we can write under $\P^{t,\bfx}_0$ that
                $$\bal
                    \vp\left({H},B^{t,\bfx_t}\right)&=\vp\left({t},\bfx^t\right)
                            +\int_t^H \Lc^{t,\bfx}\vp(s,B^{t,\bfx_t})ds
                            +\int_t^H D\vp(s,B^{t,\bfx_t})\cdot\sigma^{t,\bfx}_s(B^{t,\bfx_t})dW^{t,\bfx}_s\\
                        &= \vp\left({t},\bfx^t\right)
                            +\int_t^H \Lc^{t,\bfx}\vp(s,B^{t,\bfx_t})+D\vp(s,B^{t,\bfx_t}) \cdot f_s\nu_sds\\
                            &\quad+\int_t^H D\vp(s,B^{t,\bfx_t})\cdot\sigma^{t,\bfx}_s(B^{t,\bfx_t})\left[dW^{t,\bfx}_s-(\sigma^{t,\bfx}_s)^{-1}(B^{t,\bfx_t})f_s\nu_sds\right].
                \eal$$
                Since $\nu\in\Uc^{t,n}$,
                we have $\P^{t,\bfx}_\nu\in\Mc^{t,\bfx,\bar n}$ so that by \reff{eq: control subsol}:
                $$\bal
                    \E^{\P^{t,\bfx}_\nu}\left[(\vp-(v^n)^{t,\bfx})\left({H},B^{t,\bfx_t}\right)\right]
                        \le -c\underline\Ec_t\left[H-t\right]+\vp\left({t},\bfx^t\right)-\E^{\P^{t,\bfx}_\nu}\left[(v^n)^{t,\bfx}\left({H},B^{t,\bfx_t}\right)\right],
                \eal$$
                and taking the infimum on $\P^{t,\bfx}_\nu\in\Mc^{t,\bfx,\bar n}$ on the left-hand side and recalling that $\vp\in\underline\Ac v^n(t,\bfx)$, this gives
                \begin{align*}
                    0&<\underline\Ec_t\left[(\vp-(v^n)^{t,\bfx})\left({H},B^{t,\bfx}\right)\right]
\le -c\underline\Ec_t\left[H-t\right]
                            -\E^{\P^{t,\bfx}_\nu}\left[(v^n)^{t,\bfx}\left({H},B^{t,\bfx_t}\right)-(v^n)^{t,\bfx}\left({t},\bfx^t\right)\right],
                \end{align*}
                and finally, by the dynamic programming principle \eqref{eq: ppde for app vf dpp}, taking the infimum over $\nu\in\Uc^{t,n}$ on the right--hand side gives
                $$
                    0<-c\underline\Ec_t\left[H-t\right],
                $$
                which is a contradiction since by Lemma \ref{lemma:utile} the right--hand side is negative.

\vspace{0.5em}
            \step3{We show that $v^n$ is a viscosity super--solution to {\rm PPDE} \reff{eq: penalized PPDE}.}

\vspace{0.5em}
                Assume to the contrary that there $(t,\bfx;\vp)\in[0,T]\x\mathcal C^0\x\overline\Ac v^n(t,\bfx)$ such that for some $c>0$
                $$
                    -\Lc^{t,\bfx}\vp(t,\bfx^t)-n\rho(f_t^\top D\vp(t,\bfx^t))\le-3c<0.
                $$
                Observe again that for each $s\in[t,T]$
                $$
                    n\rho(f_s^\top D\vp(s,B^{t,\bfx}))=\sup_{u\in[0,n]^d}u\cdot (f_s^\top D\vp(s,B^{t,\bfx})),
                $$
                so that there is $u^\ast_n\in[0,n]^d$ such that
                $$
                    -\Lc^{t,\bfx}\vp(t,\bfx^t)-u^\ast_n\cdot (f_t^\top D\vp(t,\bfx^t))\le-2c.
                $$
                Without loss of generality, we may reduce $H$ in the definition of $\vp\in\overline\Ac v^n(t,\bfx)$
                so that by continuity, we obtain
                $$
                    -\Lc^{t,\bfx}\vp(t,\bfx^t)-u^\ast_n\cdot (f_s^\top D\vp(s,B^{t,\bfx_t}))\le-c,\ \mbox{on}\ [t,H],\ \mathbb P_0^{t,\bfx}-a.s.
                $$
                Set the constant control of $\Uc^{t,n}$: $\nu^n\equiv u^\ast_n$.
                Hence we have
                \beq\label{eq: control supersol}
                    -\Lc^{t,\bfx}\vp(s,B^{t,\bfx_t})-\nu^n_s\cdot f_s^\top D\vp(s,B^{t,\bfx_t})\le -c,\ \mbox{on}\ [t,H],\ \mathbb P_0^{t,\bfx}-a.s.
                \eeq
                Since $\vp$ is smooth on $[t,H]$, we have, $\P^{t,\bfx}_0-a.s.,$
                $$\bal
                    \vp\left({H},B^{t,\bfx}\right)&=\vp\left({t},\bfx^t\right)
                            +\int_t^H \Lc^{t,\bfx}\vp(s,B^{t,\bfx_t})ds
                            +\int_t^H D\vp(s,B^{t,\bfx_t})\cdot\sigma^{t,\bfx}_s(B^{t,\bfx_t})dW^{t,\bfx}_s\\
                        &= \vp\left({t},\bfx^t\right)
                            +\int_t^H \Lc^{t,\bfx}\vp(s,B^{t,\bfx_t})+D\vp(s,B^{t,\bfx_t})\cdot f_s\nu^n_sds\\
                            &\quad+\int_t^H D\vp(s,B^{t,\bfx_t})\cdot\sigma^{t,\bfx}_s(B^{t,\bfx_t})\left[dW^{t,\bfx}_s-(\sigma^{t,\bfx}_s)^{-1}(B^{t,\bfx_t})f_s\nu^n_sds\right].
                \eal$$
                By \reff{eq: control supersol}, this gives
                $$\bal
                    \E^{\P^{t,\bfx}_{\nu^n}}\left[(\vp-(v^n)^{t,\bfx})\left({H},B^{t,\bfx_t}\right)\right]
                        \ge c\E^{\P^{t,\bfx}_{\nu^n}}\left[H-t\right]
                            +\vp\left({t},\bfx^t\right)-\E^{\P^{t,\bfx}_{\nu^n}}\left[(v^n)^{t,\bfx}\left({H},B^{t,\bfx_t}\right)\right].
                \eal$$
                Since $\vp\in\overline\Ac v^n(t,\bfx)$, we have the equality $\vp({t},\bfx^t)=v^n({t},\bfx)$.
                Moreover,
                the fact that $\nu^n\in\Uc^{t,n}$ enables us to use the DPP \eqref{eq: ppde for app vf dpp} to have
                $$
                    \E^{\P^{t,\bfx}_{\nu^n}}\left[(\vp-(v^n)^{t,\bfx})\left({H},B^{t,\bfx_t}\right)\right]
                        \ge c\E^{\P^{t,\bfx}_{\nu^n}}\left[H-t\right]>0.
                $$
                Using (again) the fact that $\nu^n\in\Uc^{t,n}$ implies that $\P^{t,\bfx}_{\nu^n}\in\Mc^{t,\bfx,n}$, this contradicts
                $\vp\in\overline\Ac v^n(t,\bfx)$.

\vspace{0.5em}
\hspace{2em}$(ii)$ The case of $u^n$.

\vspace{0.5em}
            We proceed along the lines of \cite[Proof of Proposition 4.4]{ekren2016viscosity}
            and split the proof into two steps.
            Fix $(t,\bfx)\in[0,T]\x\bfL{}$ for the remainder of the proof.

\vspace{0.5em}
            \step1{We show that $u^{n}\in C^0([0,T]\times\bfL{}^{0,\bfx_0})$ and satisfies the following dynamic programming principle,
                for any $\tau\in\Tc^{t,\bfx}$:
                \beq\bal\label{eq: ppde for penalized bsde dpp}
                     \Yw^{t,\bfx,n}=
                             (u^{n})^{t,\bfx}(\tau,B^{t,\bfx_t})
                                +\int_\cdot^\tau n\rho\left[ f_s^\top ((\sigma^{t,\bfx}_s)^\top)^{-1}\left(B^{t,\bfx}\right)\Zw^{t,\bfx,n}_s \right]ds-\int_\cdot^\tau \Zw^{t,\bfx,n}_s \cdot dW^{t,\bfx}_s,
                                \;\P^{t,\bfx}_0-a.s.
                \eal\eeq}First of all, since $n\rho^\eps$ is Lipschitz--continuous and null at $0$ since Assumption \ref{assump:U} holds,
                        by standard stability results on BSDEs (see e.g. \cite{el1997backward}), for any $n\ge1$, there is a constant $C_{n}$ (which may vary from line to line) such that for all
                        $(t,\bfx,\bfx')\in[0,T]\x(\mathcal C^0)^2$
                        \begin{align}
                            &
                                \E^{\P^{t}_0}
                                \bigg[
                                    \sup_{t\le s\le T}\abs{\Yc^{t,\bfx,n}_s}^2+\int_t^T\norm{\Zc^{t,\bfx,n}_s}^2ds
                                \bigg]\le C_{n}\left( 1+\norm{\bfx}_{\infty,t}^{2(r+1)} \right),
                                    \label{eq: estimate penalized BSDE growth}&
                            \\
                            &                                \E^{\P^{t}_0}
                                \bigg[
                                    \sup_{t\le s\le T}\big|\Yc^{t,\bfx,n}_s-\Yc^{t,\bfx',n}_s\big|^2
                                    +\int_t^T\big\|\Zc^{t,\bfx,n}_s-\Zc^{t,\bfx',n}_s\big\|^2ds
                                \bigg]
                                \le C_n\big\|\bfx-\bfx'\big\|_{\infty,t}^2\big(1+\No{\bfx}_{\infty,t}^{2r}+\|\bfx'\|_{\infty,t}^{2r}\Big).
                                    \label{eq: estimate penalized BSDE Delta}&
                        \end{align}
                        In particular, this gives the following regularity
                        \be
                            &\displaystyle
                                \abs{u^{n}(t,\bfx)}\le C_{n}\left( 1+\norm{\bfx}_{\infty,t}^{r+1} \right)
                                    \label{eq: stability u deterministe growth}&\\
                            &\displaystyle
                                \mbox{and}\quad
                                \abs{u^{n}(t,\bfx)-u^{n}(t,\bfx')}\le C_n \norm{\bfx-\bfx'}_{\infty,t}\left(1+\No{\bfx}_{\infty,t}^{r}+\No{\bfx'}_{\infty,t}^{r}\right).                                    \label{eq: stability u deterministe Delta}&
                        \ee
                        By standard arguments in the BSDE theory (this would be simply the tower property for conditional expectations if $\rho$ were equal to $0$, and the result can easily be generalized using the fact that solutions to BSDEs with Lipschitz drivers can be obtained via Picard iterations), we have the following dynamic programming principle,
                        for any $t<t'\le T$
                        \begin{equation}\label{eq: ppde for penalized bsde dpp deterministe}
                            \Yw^{t,\bfx,n}=
                             (u^{n})^{t,\bfx}(t',B^{t,\bfx_t})
                                +\int_\cdot^{t'} n\rho\left[ f_s^\top((\sigma^{t,\bfx}_s)^\top)^{-1}\left(B^{t,\bfx_t}\right)\Zw^{t,\bfx,n}_s \right]ds-\int_\cdot^{t'} \Zw^{t,\bfx,n}_s \cdot dW^{t,\bfx}_s,
                                \;\P^{t,\bfx}_0-a.s.
                        \end{equation}
                        In particular $\Yw^{t,\bfx,n}_s=(u^{n})^{t,\bfx}(s,B^{t,\bfx_t})$ for any $s\in[t,T]$.
                        It then follows                         
                        \begin{align}
                            \displaystyle
                            &|u^n(t,\bfx)-u^n(t',\bfx)|
                                = \Big|
                                        \E^{\P^{t,\bfx}_0}\Big[
                                            \Yw^{t,\bfx,n}_t
                                            -\Yw^{t,\bfx,n}_{t'}
                                            +(u^{n})^{t,\bfx}(t',B^{t,\bfx_t})
                                            -u^n(t',\bfx)
                                        \Big]\Big|\nonumber\\
                            \displaystyle
                               \nonumber &\le
                                        \E^{\P^{t,\bfx}_0}\bigg[
                                            \int_t^{t'} \big| n\rho\big(f_s^\top((\sigma^{t,\bfx}_s)^\top)^{-1}(B^{t,\bfx_t})\Zw^{t,\bfx,n}_s \big)\big|ds\bigg]\nonumber     +\E^{\P^{t,\bfx}_0}\left[\abs{(u^{n})^{t,\bfx}(t',B^{t,\bfx_t})
                                            -u^n(t',\bfx)}
                                        \right]\nonumber\\
                            \displaystyle
                           \nonumber     &\le  \E^{\P^{t,\bfx}_0}\bigg[
                                            \int_t^{t'} \big| n\rho\big(f_s^\top((\sigma^{t,\bfx}_s)^\top)^{-1}(B^{t,\bfx_t})\Zw^{t,\bfx,n}_s \big)\big|ds\bigg]\\
                                            &\hspace{0.9em}+C_n\bigg(\underset{t\leq s\leq t'}{\sup}\No{\bfx_s-\bfx_t}+\mathbb E^{\mathbb P_0^{t,\bfx}}\bigg[\underset{t\leq s\leq t'}{\sup}\big\|B^{t,\bfx_t}_s-\bfx_t\big\|^2\bigg]^{\frac12}\bigg)\bigg(1+\No{\bfx}_{\infty,t'}^r+\mathbb E^{\mathbb P_0^{t,\bfx}}\bigg[\underset{t\leq s\leq t'}{\sup}\big\|B^{t,\bfx_t}_s\big\|^2\bigg]^{\frac12}\Bigg)\nonumber\\
                            \displaystyle
                                &\le \E^{\P^{t,\bfx}_0}\bigg[
                                            \int_t^{t'} \abs{n\rho\big( f_s^\top((\sigma^{t,\bfx}_s)^\top)^{-1}(B^{t,\bfx_t})\Zw^{t,\bfx,n}_s \big)}ds
                                 +C_nd_\infty((t,\bfx);(t',\bfx))\left(1+\No{\bfx}_{\infty,t'}^r\right),
                                        \label{eq: for DPP regularity}
                        \end{align}
                        where the last line follows from \reff{eq: estimate SDE Delta initial}.
                        Observe from \reff{eq: estimate penalized BSDE growth} combined with Assumption \ref{ass: mu sigma lip} and the Lipschitz-continuity of $\rho^\eps$
                        that
                        \begin{align*}
                            \E^{\P^{t,\bfx}_0}\left[
                                            \int_t^{t'} \abs{n\rho\left( f_s^\top((\sigma^{t,\bfx}_s)^\top)^{-1}\left(B^{t,\bfx_t}\right)\Zw^{t,\bfx,n}_s \right)}ds
                                        \right]&= \E^{\P^{t}_0}\left[
                                            \int_t^{t'} \abs{n\rho\left( f_s^\top((\sigma^{t,\bfx}_s)^\top)^{-1}\left(B^{t}\right)\Zc^{t,\bfx,n}_s \right)}ds
                                        \right]\\
                                &\le C_n \left( 1+\No{ \bfx}_{\infty,t'}^{r+1} \right)(t'-t)^{\frac12}.
                        \end{align*}
                        Plugging the latter into \reff{eq: for DPP regularity} and since $\sqrt{t'-t}\le d_\infty((t,\bfx);(t',\bfx))$, this gives finally
                        $$
                            \abs{u^n(t,\bfx)-u^n(t',\bfx)}
                            \le C_n \left( 1+\No{ \bfx}_{\infty,t'}^{r+1} \right) d_\infty((t,\bfx);(t',\bfx)).
                        $$
                        Finally, the regularity in time that we just proved, allows us to classically extend the dynamic programming principle in \reff{eq: ppde for penalized bsde dpp deterministe}
                        to stopping times, giving \reff{eq: ppde for penalized bsde dpp} (the result is clear for stopping times taking finitely many values, and the general result follows by the usual approximation of stopping times by decreasing sequences of stopping times with finitely many values).
                        
                        \vspace{0.5em}
                    \step2{We conclude the proof.}

                        Without loss of generality,
                        we prove only the viscosity sub--solution,
                        the super--solution being obtained similarly.
                        Assume to the contrary that there is $(t,\bfx;\vp)\in[0,T]\x\mathcal C^0\x\underline\Ac u^n(t,\bfx)$
                        such that
                        $$
                            2c:=-\Lc^{t,\bfx}\vp(t,\bfx^t)
                                -n\rho\left(f_t^\top D\vp(t,\bfx^t)\right)>0.
                        $$
                        Let $H$ be the hitting time corresponding to the definition of $\vp\in\underline\Ac u^n(t,\bfx)$.
                        By continuity of $\vp$ and $\rho$,  
                        reducing $H$ if necessary, we deduce that
                        $$
                            -\Lc^{t,\bfx}\vp(s,B^{t,\bfx_t})
                                -n\rho\left(f_sD\vp(s,B^{t,\bfx_t})\right)\ge c>0
                               ,\ s\in[t,H].
                        $$
                        By
                        the DPP \reff{eq: ppde for penalized bsde dpp} and the smoothness of $\vp$, we have under $\P_0^{t,\bfx}$,
                        we have
                        \begin{align*}
                          &(\vp-(u^{n})^{t,\bfx})_H(\cdot,B^{t,\bfx_t})
                                -(\vp-(u^{n})^{t,\bfx})_t(\cdot,\bfx^t)\\
                            &=\int_t^H \Lc^{t,\bfx}\vp(s,B^{t,\bfx_t})ds
                                +\int_t^H \sigma^{t,\bfx}_s\left(B^{t,\bfx_t}\right) D\vp(s,B^{t,\bfx_t})\cdot dW^{t,\bfx}_s\\
                            &\quad
                                +\int_t^H n\rho\left(f_s^\top \left(( \sigma^{t,\bfx}_s)^\top \right)^{-1}\left(B^{t,\bfx_t}\right)\Zw^{t,\bfx,n}_s\right)ds
                                    - \int_t^H \Zw^{t,\bfx,n}_s\cdot dW^{t,\bfx}_s\\
                            &\le -\int_t^H
                                \left[
                                    c+n\left(
                                            \rho\left(f_s^\top D\vp(s,B^{t,\bfx_t})\right)
                                            -\rho\left(f_s^\top\left( (\sigma^{t,\bfx}_s)^\top\right)^{-1}\left(B^{t,\bfx_t}\right)\Zw^{t,\bfx,n}_s\right)
                                        \right)
                                \right]ds\\
                                &\quad+\int_t^H \sigma^{t,\bfx}_s\left(B^{t,\bfx_t}\right)
                                    \left(
                                         D\vp(s,B^{t,\bfx_t})
                                        -\left(\sigma^{t,\bfx}_s\right)^{-1}\left(B^{t,\bfx_t}\right) \Zw^{t,\bfx,n}_s
                                    \right)\cdot
                                        dW^{t,\bfx}_s\\
                            &= -c(H-t)
                                +\int_t^H\alpha^n_s\cdot
                                    \left(
                                        D\vp(s,B^{t,\bfx_t})
                                        -\left( \sigma^{t,\bfx}_s\right)^{-1}\left(B^{t,\bfx_t}\right)\Zw^{t,\bfx,n}_s
                                    \right)ds\\
                                &\quad+\int_t^H \sigma^{t,\bfx}_s\left(B^{t,\bfx_t}\right)
                                    \left(
                                         D\vp(s,B^{t,\bfx_t})
                                        -\left(\sigma^{t,\bfx}_s\right)^{-1}\left(B^{t,\bfx_t}\right) \Zw^{t,\bfx,n}_s
                                    \right)\cdot
                                        dW^{t,\bfx}_s\\
                            &= c(t-H)
                                +\int_t^H
                                    \left(
                                        D\vp(s,B^{t,\bfx_t})
                                        -\left( \sigma^{t,\bfx}_s\right)^{-1}\left(B^{t,\bfx_t}\right)\Zw^{t,\bfx,n}_s
                                    \right)\cdot
                                    \left(
                                        \sigma^{t,\bfx}_s\left(B^{t,\bfx_t}\right) dW^{t,\bfx}_s
                                        +\alpha^n_sds
                                    \right),
                        \end{align*}
                        where $|\alpha^n|\le n||f||$
                        $ds\otimes d\P^{t,\bfx}_0-a.e.$
                        By Girsanov's Theorem, we then have that there is $\Q\in\Mc^{t,\bfx,\bar n}$ such that
                        $\int\sigma^{t,\bfx}_s(B^{t,\bfx_t})dW^{t,\bfx}_s+\alpha^n_sds$ is a $\Q-$Brownian motion.
                        The above inequality holds then $\Q^n-a.s.$
                        so that
                        $$\bal
                                (\vp-(u^{n})^{t,\bfx})_t(\cdot,\bfx^t)&\ge
                                \E^{\Q^n}\left[(\vp-(u^{n})^{t,\bfx})_H(\cdot,B^{t,\bfx_t})+c(H-t)\right]>\E^{\Q^n}\left[(\vp-(u^{n})^{t,\bfx})_H(\cdot,B^{t,\bfx_t})\right],
                        \eal$$
                        which is in contradiction with the definition of $\vp\in\underline\Ac u^{n}(t,\bfx)$.
                \ep
                
\section{Technical proofs for Section \ref{Sec5}}

\begin{lemma}\label{lem.operator}
Let Assumption \ref{ass: mu sigma lip2} hold. Fix $n\geq 1$. For every $(t,s,\bfx,\tilde\bfx,\lambda,u)\in[0,T]\times[t,T]\times\mathcal C^0\times\mathcal C^t\times\R_+^*\times\R$, for every $k=0,\dots,n-1$, and for every Borel map $f:\R^d\longrightarrow\R$ with polynomial growth, define the following operators
\begin{align*}
&Q^{k+1}_{t,\bfx,\lambda}(f)(\tilde\bfx,u):=\E^{\P_0^t}\Big[f\Big.\Big(\tilde\bfx_{t^{t,n}_k}+\lambda\Big(\mu_{t^{t,n}_k}^{t,\bfx}(\tilde\bfx)+f_{t^{t,n}_k}\nu_{t^{t,n}_k}\Big)+\Big(uM_{t^{t,n}_k}+\sigma^{t,\bfx}_{t^{t,n}_k}(\tilde\bfx)\Big)\Big(B^t_{t^{t,n}_{k+1}}-B^t_{t^{t,n}_k}\Big)\Big)\Big|\mathcal F^{t}_{t^{t,n}_k}\Big].
\end{align*}
If $f$ is concave and s.t. for every $(t,s,\bfx,\tilde\bfx,\tilde\bfy,\alpha,\beta,\gamma,\eta)\in[0,T]\times[t,T]\times\mathcal C^0\times\mathcal C^t\times\mathcal C^t\times\R^d\times\R_+^*\times\R^d\times[0,1]$, 
\begin{align}\label{eq:f}
\nonumber&f\left(\alpha+\beta\mu^{t,x}_s(\eta\tilde\bfx+(1-\eta)\tilde\bfy)+\sigma_s^{t,\bfx}(\eta\tilde\bfx+(1-\eta)\tilde\bfy)\gamma\right)
\\&\geq  f\left(\alpha+\beta(\eta\mu^{t,x}_s(\tilde\bfx)+(1-\eta)\mu^{t,x}_s(\tilde\bfy))+(\eta\sigma_s^{t,\bfx}(\tilde\bfx)+(1-\eta)\sigma^{t,x}_s(\tilde\bfy))\gamma\right),
\end{align}
then the map $(\tilde\bfx,u)\longmapsto Q^{k+1}_{t,\bfx,\lambda}(f)(\tilde\bfx,u)$ is concave,
and the map $u\longmapsto Q^{k+1}_{t,\bfx,\lambda}(f)(\tilde\bfx,u)$ is non--decreasing on $\R_-$.
\end{lemma}
  \proof
  The fact that the operators $Q^{k+1}_{t,\bfx,\gamma}$ are well-defined is clear from the polynomial growth of $f$, the linear growth of $\mu$ and $\sigma$, and the fact that $B^t$ has moments of any order under $\P_0^t$. Then, the concavity of $Q^{k+1}_{t,\bfx,\gamma}(f)$ is an immediate corollary of the concavity assumptions on $f$, as well as \reff{eq:f}. 

\vspace{0.5em}
\noindent Then, since $f$ has polynomial growth, Feynman--Kac's formula implies that 
$$ Q^{k+1}_{t,\bfx,\lambda}(f)(\tilde\bfx,u)=v(t^{t,n}_{k},B^t_{t^{t,n}_k}),$$ 
where $v:[t^{t,n}_k,t^{t,n}_{k+1}]\times\R^d\longrightarrow \R$ is the unique viscosity solution of the PDE
$$\begin{cases}
\displaystyle -v_s-\frac12{\rm Tr}\Big[\Big(uM_{t^{t,n}_k}+\sigma^{t,\bfx}_{t^{t,n}_k}(\tilde\bfx)\Big)\Big(uM_{t^{t,n}_k}+(\sigma^{t,\bfx}_{t^{t,n}_k})^\top(\tilde\bfx)\Big)v_{xx}\Big]=0,\ \text{on }[t^{t,n}_k,t^{t,n}_{k+1})\times\R^d,\\[0.5em]
\displaystyle v\left(t^{t,n}_{k+1},x\right)=f\Big(\tilde\bfx_{t^{t,n}_k}+\lambda\Big(\mu_{t^{t,n}_k}^{t,\bfx}(\tilde\bfx)+f_{t^{t,n}_k}\nu_{t^{t,n}_k}\Big)+\Big(uM_{t^{t,n}_k}+\sigma^{t,\bfx}_{t^{t,n}_k}\Big)x\Big),\ x\in\R^d.
\end{cases}$$
 This linear PDE classically satisfies a comparison theorem, and $v$ is concave in $x$ because of the concavity of $f$. Moreover, the diffusion part of the PDE rewrites, as a quadratic functional of $u$
 $$u^2\Tr{M_{t^{t,n}_k}^2v_{xx}}+u\Tr{\left(M_{t^{t,n}_k}(\sigma^{t,\bfx}_{t^{t,n}_k})^\top(\tilde\bfx)+\sigma^{t,\bfx}_{t^{t,n}_k}(\tilde\bfx)M_{t^{t,n}_k}\right)v_{xx}}+\Tr{\sigma^{t,\bfx}_{t^{t,n}_k}(\tilde\bfx)(\sigma^{t,\bfx}_{t^{t,n}_k})^\top(\tilde\bfx)}.$$
 Since $M^2_t$ is symmetric positive, $M_t\sigma^\top_t(\bfx)+\sigma_t(\bfx)M_t$ is symmetric negative and $v_{xx}$ is symmetric negative as well, the above is actually non-decreasing for $u\in \R_-$. The same then holds for $Q^{k+1}_{t,\bfx,\lambda}(f)(\tilde\bfx,u)$ by comparison.
      \ep
   \begin{proposition}\label{prop:dege}
Let Assumptions \ref{assump:U} and \ref{ass: mu sigma lip2} hold and fix some $q\geq p> 0$. For every $n\in\mathbb N\backslash\{0\}$, $\bfx\in\mathcal C^0$, let us define recursively $(X_k^{t,\bfx,n})_{0\leq k\leq n}$ and $(Y_k^{t,\bfx,n})_{0\leq k\leq n}$ by $X_0^{t,\bfx,n}=Y_0^{t,\bfx,n}={\bf x}_t$, and for $0\leq k\leq n-1$
   \begin{align*}
   X_{k+1}^{t,\bfx,n}=&\ X_k^{t,\bfx,n}+\Big(\mu^{t,\bfx}_{t_k^{t,n}}\big(\overline X_k^{t,\bfx,n}\big)+f_{t_k^{t,n}}\nu_{t^{t,n}_k}\Big)\big(t^{t,n}_{k+1}-t^{t,n}_k\big)+\Big(\eta^{p,t,\bfx}_{t^{t,n}_k}\big(\overline X^{t,\bfx,n}_k\big)M_{t^{t,n}_k}+\sigma^{t,\bfx}_{t^{t,n}_k}\big(\overline X^{t,\bfx,n}_k\big)\Big)\big(B^t_{t^{t,n}_{k+1}}-B^t_{t^{t,n}_k}\big),\\
   Y_{k+1}^{t,\bfx,n}=&\ Y_k^{t,\bfx,n}+\Big(\mu^{t,\bfx}_{t_k^{t,n}}\big(\overline Y_k^{t,\bfx,n}\big)+f_{t_k^{t,n}}\nu_{t^{t,n}_k}\Big)\big(t^{t,n}_{k+1}-t^{t,n}_k\big)+\Big(\eta^{q,t,\bfx}_{t^{t,n}_k}\big(\overline Y^{t,\bfx,n}_k\big)M_{t^{t,n}_k}+\sigma^{t,\bfx}_{t^{t,n}_k}\big(\overline Y^{t,\bfx,n}_k\big)\Big)\big(B^t_{t^{t,n}_{k+1}}-B^t_{t^{t,n}_k}\big),
   \end{align*}
  where $(\overline X^{t,\bfx, n}_k)_{0\leq k\leq n}$ and $(\overline Y^{t,\bfx,n}_k)_{0\leq k\leq n}$ are defined as the following piecewise linear interpolations 
   $$\overline X^{t,\bfx, n}_k:=i_k(X^{t,\bfx,n}_{0:k}),\ \overline Y^{t,\bfx, n}_k:=i_k(Y^{t,\bfx,n}_{0:k}).$$
   
Then, we have
$$\E^{\P_0^t}\left[U^{t,\bfx}\left(i_n\left(X_{0:n}^{t,n,\bfx}\right)\right)\right]\leq \E^{\P_0^t}\left[U^{t,\bfx}\left(i_n\left(Y_{0:n}^{t,n,\bfx}\right)\right)\right].$$
   \end{proposition}
   
   \proof
  Let $\Delta^{t,n}:=t^{t,n}_{k+1}-t^{t,n}_k=(T-t)/n$, and consider the following martingales, for $0\leq k\leq n$,
   $$M_k:=\E^{\mathbb P_0^t}\left[\left.U^{t,\bfx}\left(i_n\left(X^{t,n,\bfx}_{0:n}\right)\right)\right|\mathcal F_{t^{t,n}_k}^t\right],\ N_k:=\E^{\mathbb P_0^t}\left[\left.U^{t,\bfx}\left(i_n\left(Y^{t,n,\bfx}_{0:n}\right)\right)\right|\mathcal F_{t^{t,n}_k}^t\right],$$
   which are well-defined, since $U$ has polynomial growth and we know from Lemma \ref{lem.estimates} that $X^{t,n,\bfx}_{0:n}$ and $Y^{t,n,\bfx}_{0:n}$ have moments of any order. For any $k=0,\dots,n$, we also define the following sequences of functions from $\R^{k+1}$ to $\R$, for $k=0,\dots,n-1$, by backward induction, for any $x_{0:k}\in(\R^d)^{k+1}$
   \begin{align*}
  \Phi_n:=&\ U^{t,\bfx}\circ i_n,\ \Phi_k(x_{0:k}):=Q^{k+1}_{t,\bfx,\Delta^{t,n}}(\Phi_{k+1}(x_{0:k},\cdot))\Big(i_k(x_{0:k}),\eta^p_{t^{t,n}_k}(i_k(x_{0:k}))\Big),\\
   \Psi_n:=&\ U^{t,\bfx}\circ i_n,\ \Psi_k(x_{0:k}):=Q^{k+1}_{t,\bfx,\Delta^{t,n}}(\Psi_{k+1}(x_{0:k},\cdot))\Big(i_k(x_{0:k}),\eta^q_{t^{t,n}_k}(i_k(x_{0:k}))\Big).
      \end{align*}
      It is immediate by definition of $X^{t,\bfx,n}$ and $Y^{t,\bfx,n}$ that we have for every $0\leq k\leq n$
      $$M_k=\Phi_k(X^{t,\bfx,n}_{0:k})\ \text{and}\ N_k=\Psi_k(Y^{t,\bfx,n}_{0:k}).$$
      Let us now show that the maps $\Phi_k$ and $\Psi_k$ are concave for every $k=0,\dots,n$, and that they verify that for any $0\leq i\leq k-1$, $(x_{0:n},\tilde x_{0:n},\eta)\in(\R^d)^{n+1}\times(\R^d)^{n+1}\times[0,1]$, for any $\{(\alpha_{m,l},\beta_{m,l},\gamma_{m,l})\in\R^d\times\times \R^*_+\times\R^d,\ i\leq m\leq k-1,\ 0\leq l\leq k-i-1\}$, we have for $\varphi=\Phi,\Psi$
\begin{align}\label{eq:eq}
\nonumber&\varphi_k\Big(x_{0:i},\alpha_{i,0} +\beta_{i,0}\mu^{t,\bfx}_{t^{t,n}_{i}}(i_{i}(\eta x_{0:i}+(1-\eta)\tilde x_{0:i})+\sigma^{t,\bfx}_{t^{t,n}_{i}}(i_{i}(\eta x_{0:i}+(1-\eta)\tilde x_{0:i})\gamma_{i,0},\\
\nonumber &\Sum_{j=i}^{i+1}\Big(\alpha_{j,1}+\beta_{j,1}\mu^{t,\bfx}_{t^{t,n}_{j}}(i_{j}(\eta \hat x_{0:j}+(1-\eta)\hat{\tilde x}_{0:j})+\sigma^{t,\bfx}_{t^{t,n}_{j}}(i_{j}(\eta \hat x_{0:j}+(1-\eta)\hat{\tilde x}_{0:j})\gamma_{j,1}\Big)\Big),\dots, \\
\nonumber &\Sum_{j=i}^{k-1}\Big(\alpha_{j,k-i-1}+\beta_{j,k-i-1}\mu^{t,\bfx}_{t^{t,n}_{j}}(i_{j}(\eta \hat x_{0:j}+(1-\eta)\hat{\tilde x}_{0:j})+\sigma^{t,\bfx}_{t^{t,n}_{j}}(i_{j}(\eta \hat x_{0:j}+(1-\eta)\hat{\tilde x}_{0:j})\gamma_{j,k-i-1}\Big)\Big)\\
\nonumber&\geq      \varphi_k\Big(x_{0:i},\alpha_{i,0} +\beta_{i,0}\Big(\eta\mu^{t,\bfx}_{t^{t,n}_{i}}(i_{i}(x_{0:i}))+(1-\eta)\mu^{t,\bfx}_{t^{t,n}_{i}}(i_{i}(\tilde x_{0:i}))\Big)+\Big(\eta\sigma^{t,\bfx}_{t^{t,n}_{i}}(i_{i}(x_{0:i}))+(1-\eta)\sigma^{t,\bfx}_{t^{t,n}_{i}}(i_{i}(\tilde x_{0:i}))\Big)\gamma_{i,0}\Big),\\
\nonumber&\Sum_{j=i}^{i+1}\Big(\alpha_{j,1} +\beta_{j,1}\Big(\eta\mu^{t,\bfx}_{t^{t,n}_{j}}(i_{j}(\hat x_{0:j}))+(1-\eta)\mu^{t,\bfx}_{t^{t,n}_{j}}(i_{j}(\hat{\tilde x}_{0:j}))\Big)+\Big(\eta\sigma^{t,\bfx}_{t^{t,n}_{j}}(i_{j}(\hat x_{0:j}))+(1-\eta)\sigma^{t,\bfx}_{t^{t,n}_{j}}(i_{j}(\hat{\tilde x}_{0:j}))\Big)\gamma_{j,1}\Big),\dots,\\
&\Sum_{j=i}^{k-1}\Big(\alpha_{j,k-i-1} +\beta_{j,k-i-1}\Big(\eta\mu^{t,\bfx}_{t^{t,n}_{j}}(i_{j}(\hat x_{0:j}))+(1-\eta)\mu^{t,\bfx}_{t^{t,n}_{j}}(i_{j}(\hat{\tilde x}_{0:j}))\Big)+\Big(\eta\sigma^{t,\bfx}_{t^{t,n}_{j}}(i_{j}(\hat x_{0:j}))+(1-\eta)\sigma^{t,\bfx}_{t^{t,n}_{j}}(i_{j}(\hat{\tilde x}_{0:j}))\Big)\gamma_{j,k-i-1}\Big),
\end{align}
      where $\hat x$ and $\hat{\tilde x}$ are defined recursively, for $w:=x,\tilde x$, by
   $$   \begin{cases}
    \displaystyle  \hat w_l:=w_l,\text{ $0\leq l\leq i,$ }\\
      \displaystyle\hat{w}_{l+1}:=\Sum_{j=i}^{l}\Big(\alpha_{j,l-i}+\beta_{j,l-i}\mu^{t,\bfx}_{t^{t,n}_{j}}(i_{j}( \hat w_{0:j}))+\sigma^{t,\bfx}_{t^{t,n}_{j}}(i_{j}( \hat x_{0:j}))\gamma_{j,l-i}\Big)\Big),\ i\leq l\leq k-1.
      \end{cases}$$
      We only prove the result for $\Phi_k$, the other one being exactly similar. We argue by backward induction. When $k=n$, the result is obvious since $U$ is concave and Assumption \ref{ass: mu sigma lip2}(vi) holds. Let us assume that the properties holds for $\Phi_{k+1}$ for some $k\leq n-1$. Then, let us now show that the map $x_{0:k}\longmapsto Q^{k+1}_{t,\bfx,\Delta^{t,n}}(\Phi_{k+1}(x_{0:k},\cdot))\left(i_k(x_{0:k}),u\right)$ is concave for any $u\in\R$. We actually have
      \begin{align*}
     & Q^{k+1}_{t,\bfx,\Delta^{t,n}}(\Phi_{k+1}(x_{0:k},\cdot))\left(i_k(x_{0:k}),u\right)\\
      &=\E^{\P_0^t}\Big[\Phi_{k+1}\Big(x_{0:k},x_k+\Delta^{t,n}\Big(\mu_{t^{t,n}_k}^{t,\bfx}(i_k(x_{0:k}))+f_{t^{t,n}_k}\nu_{t^{t,n}_k}\Big)+\Big(uM_{t^{t,n}_k}+\sigma^{t,\bfx}_{t^{t,n}_k}(i_k(x_{0:k}))\Big)\Big(B^t_{t^{t,n}_{k+1}}-B^t_{t^{t,n}_k}\Big)\Big|\mathcal F^t_{t^{t,n}_k}\Big].
      \end{align*}
      Therefore, the concavity is immediate from the induction hypothesis on $\Phi_{k+1}$ (both the concavity and Inequality \reff{eq:eq}). Now, we know that $\eta_t^p(\cdot)$ is concave and non--positive, and, from Lemma \ref{lem.operator}, that the map $u\longmapsto Q^{k+1}_{t,\bfx,\Delta^{t,n}}(\Phi_{k+1}(x_{0:k},\cdot))\left(i_k(x_{0:k}),u\right)$ is non--decreasing on $\R_-$. This therefore proves the concavity of $\Phi_k$. Moreover, $\Phi_k$ inherits \reff{eq:eq} directly from $\Phi_{k+1}$ by its definition as an expectation of $\Phi_{k+1}$. 
      
      \vspace{0.5em}
      \noindent Finally, let us prove, again by backward induction, that for every $k=0,\dots,n$, $\Phi_k\leq \Psi_k$. The result is obvious by definition for $k=n$. Assume now that for some $k\leq n-1$, we have $\Phi_{k+1}\leq \Psi_{k+1}$. Then, for any $x_{0:k}\in\R^{k+1}$
      \begin{align*}
      \Phi_k(x_{0:k})&=Q^{k+1}_{t,\bfx,\Delta^{t,n}}(\Phi_{k+1}(x_{0:k},\cdot))\left(i_k(x_{0:k}),\eta^p_{t^{t,n}_k}(i_k(x_{0:k}))\right)\\
      &\leq Q^{k+1}_{t,\bfx,\Delta^{t,n}}(\Phi_{k+1}(x_{0:k},\cdot))\left(i_k(x_{0:k}),\eta^q_{t^{t,n}_k}(i_k(x_{0:k}))\right)\\
      &\leq Q^{k+1}_{t,\bfx,\Delta^{t,n}}(\Psi_{k+1}(x_{0:k},\cdot))\left(i_k(x_{0:k}),\eta^q_{t^{t,n}_k}(i_k(x_{0:k}))\right)=\Psi_k(x_{0:k}),
      \end{align*}
      where we have used successively the fact that $u\longmapsto Q^{k+1}_{t,\bfx,\Delta^{t,n}}(\Phi_{k+1}(x_{0:k},\cdot))\left(i_k(x_{0:k}),u\right)$ is non--decreasing on $\R_-$ (remember that $\eta^p\leq \eta^q$) and the induction hypothesis. To conclude, it suffices to take $k=0$ to obtain $\Phi_0(\bfx_t)\leq \Psi_0(\bfx_t),$ which is equivalent by the martingale property of $M$ and $N$ to 
      $$\E^{\P_0^t}\left[U^{t,\bfx}\left(i_n\left(X_{0:n}^{t,n,\bfx}\right)\right)\right]\leq \E^{\P_0^t}\left[U^{t,\bfx}\left(i_n\left(Y_{0:n}^{t,n,\bfx}\right)\right)\right].$$
   \ep
   
   We can now state the main technical result of this section.
   
   \begin{proposition}\label{prop:decreasing}
  Let Assumptions \ref{assump:U} and \ref{ass: mu sigma lip2} hold. For any $p> 0$, denote by $X^{t,\bfx,\nu,p}$ the solution to the {\rm SDE} \reff{eq:sde} with diffusion matrix $\sigma^p$ instead of $\sigma$, and let $$v^p(t,\bfx):=\underset{\nu\in\mathcal U^t}{\sup}\mathbb E^{\P_0^t}\left[U\left(\bfx\otimes_tX^{t,\bfx,\nu,p}\right)\right].$$ Then, for any $q\geq p> 0$, we have for any $(t,\bfx)\in[0,T]\times\mathcal C^0$
  $$v^p(t,\bfx)\leq v^q(t,\bfx).$$
   \end{proposition}
   \proof
   By Proposition \ref{prop:dege}, we know that if we replace $X^{t,\bfx,\nu,p}$ and $X^{t,\bfx,\nu,q}$ by their Euler scheme, then the expectation of $U$ of these Euler schemes are ordered. We can then follow exactly the arguments of the proofs of Lemma 2.2 and Theorem 2.1 in \cite{pages2014convex}, using in particular the continuity we have assumed for $U$, as well as the fact that the genuine Euler scheme for a non--Markovian SDE converges to the solution of the SDE for the uniform topology on $\mathcal C^0$, to extend this result and obtain
   $$\mathbb E^{\P_0^t}\left[U\left(\bfx\otimes_tX^{t,\bfx,\nu,p}\right)\right]\leq \mathbb E^{\P_0^t}\left[U\left(\bfx\otimes_tX^{t,\bfx,\nu,q}\right)\right],$$
   from which the result is clear.
   \ep

\section{Technical proofs for Section \ref{Sec6}}

\proof[Proof of Lemma \ref{lemma:face}]
Clearly, we have $\widehat U\geq U$, so that one inequality is trivial. Next, fix some $u\in\mathcal V^t_{\rm b}$. For some $\eps_o>0$ small enough and any $\eps\in(0,\eps_o)$, we define the following element of $\mathcal V^t_{\rm b}$
$$u^\eps_s:= u_s{\bf 1}_{[t,T-\eps]}(s)+\frac{\iota}{\eps}{\bf 1}_{[T-\eps,T]}(s),$$
where $\iota$ is any bounded and $\mathcal F^t_{T-\eps_o}-$measurable random variable. Then, we have by the tower property for expectations and the Markov property for SDEs (see for instance \cite{claisse2016pseudo}) that
\begin{align*}
&\mathbb E^{\P_0^t}\bigg[\mathbb E^{\P_0^{T-\eps}}\left[U\bigg(X_T^{T-\eps,X_{T-\eps}^{t,x,u},(u^\eps)^{T-\eps,B^t}}\bigg)-\int_{T-\eps}^T\delta\big((u^\eps)^{T-\eps,B^t}_s\big)ds\right]-\int_t^{T-\eps}\delta(u_s)ds\bigg]\\
&=\mathbb E^{\P_0^t}\left[U(X_T^{t,x,u^\eps})-\int_{t}^T\delta(u^\eps_s)ds\right],
\end{align*}
which implies that
\begin{align}
 Y_t^\bfx\geq \mathbb E^{\P_0^t}\bigg[\mathbb E^{\P_0^{T-\eps}}\left[U\bigg(X_T^{T-\eps,X^{t,x,u}_{T-\eps},(u^\eps)^{T-\eps,B^t}}\bigg)-\int_{T-\eps}^T\delta\big((u^\eps)^{T-\eps,B^t}_s\big)ds\right]-\int_t^{T-\eps}\delta(u_s)ds\bigg].\label{eq:eqeqeq}
\end{align}
Next, we claim that, at least along a subsequence, we have
\begin{align}\label{eq:eqeq}
\nonumber&\underset{\eps\rightarrow 0}{\rm lim}\ \mathbb E^{\P_0^t}\bigg[\mathbb E^{\P_0^{T-\eps}}\left[U\bigg(X_T^{T-\eps,X^{t,x,u}_{T-\eps},(u^\eps)^{T-\eps,B^t}}\bigg)-\int_{T-\eps}^T\delta\big((u^\eps)^{T-\eps,B^t}_s\big)ds\right]-\int_t^{T-\eps}\delta(u_s)ds\bigg]\\
&\hspace{4.2em}= \mathbb E^{\P_0^t}\left[U\big(X_T^{t,x,u}+f_T\iota\big)-\delta(\iota)-\int_t^T\delta(u_s)ds\right].
\end{align}
Indeed, by (an easy extension of) \eqref{eq: estimate SDE Delta}, we first have for any $p\geq 2$
\begin{align*}
\E^{\P^t_0}
                                \left[
                                   \norm{X_T^{T-\eps,X_{T-\eps}^{t,x,u},u^{T-\eps,B^t}}+f\iota-X_T^{T-\eps,X_{T-\eps}^{t,x,u},(u^\eps)^{T-\eps,B^t}}}^p
                                \right]
                                &\le C_p\E^{\P_0^t}\left[\No{f\iota+\int_{T-\eps}^Tf\left(u_s^{T-\eps,B^t}-\frac{\iota}{\eps}\right)ds}^{p}\right]\\
                                &\leq C_p\int_{T-\eps}^T\E^{\P_0^t}\left[\No{fu_s^{T-\eps,B^t}}^p\right]ds.
                                \end{align*}
                                Furthermore, by \eqref{eq: estimate SDE Delta initial}, we have $\P_0^t{\rm -a.s.}$
                                 \begin{align*}
                                \E^{\P^{T-\eps}_0}
                                \Big[
                                    \Big\|X_T^{T-\eps,X_{T-\eps}^{t,x,u},u^{T-\eps,B^t}}-X^{t,x,u}_{T-\eps}\Big\|^p
                                \Big]\le&\ C_p\eps^{\frac12}\left(1+\No{X_{T-\eps}^{t,x,u}}^p\right)+C_p\E^{\P_0^{T-\eps}}\left[\No{\int_{T-\eps}^{T}u_s^{T-\eps,B^t}ds}^{p}\right],
\end{align*}
which implies by the tower property
                                 \begin{align*}
\displaystyle
                                \E^{\P^{t}_0}
                                \Big[
                                    \Big\|X_T^{T-\eps,X_{T-\eps}^{t,x,u},u^{T-\eps,B^t}}-X^{t,x,u}_{T-\eps}\Big\|^p
                                \Big]\le&\ C_p\eps^{\frac12}\left(1+\E^{\P_0^{t}}\left[\No{X^{t,x,u}_{T-\eps}}^p\right]\right)+C_p\E^{\P_0^{t}}\left[\No{\int_{T-\eps}^{T}u_sds}^{p}\right],\ \P_0^t{\rm -a.s.}
\end{align*}
                            Hence, passing to a subsequence if necessary, and using the continuity of the paths of $X^{t,x,u}$, we deduce from the above inequalities that
                               $$\text{$ X_T^{T-\eps,X^{t,x,u},(u^\eps)^{T-\eps,B^t}}$ converges to $X^{t,x,u}_T+f\iota$, $\P_0^t-a.s.$ and in $L^p(\mathbb P_0^t)$.}$$
                                By continuity of $U$, we deduce that the following convergence holds $\P_0^t-a.s. $ and in $L^p(\mathbb P_0^t)$
                                 $$\text{$ U\bigg(X_T^{T-\eps,X_{T-\eps}^{t,x,u},(u^\eps)^{T-\eps,B^t}}\bigg)-\int_{T-\eps}^T\delta\big((u^\eps)^{T-\eps,B^t}_s\big)ds\longrightarrow U(X^{t,x,u}_T+f\iota)-\delta(\iota)$.}$$
Then, this implies by the tower property that the following convergence holds in $L^1(\mathbb P_0^t)$
 \begin{align*}
 & \E^{\P_0^{T-\eps}}\left[U\bigg(X_T^{T-\eps,X_{T-\eps}^{t,x,u},(u^\eps)^{T-\eps,B^t}}\bigg)-\int_{T-\eps}^T\delta\big((u^\eps)^{T-\eps,B^t}_s\big)ds\right]\longrightarrow U(X^{t,x,u}_T+f\iota)-\delta(\iota),
 \end{align*}
 which implies the desired claim \eqref{eq:eqeq}.

\vspace{0.5em}
Then, by \eqref{eq:eqeqeq} and \eqref{eq:eqeq} we deduce that for random variable $\iota$ which is bounded and $\Fc^t_{T-\eps_o}-$measurable, we have
\begin{equation}\label{eq:face}
Y_t^x\geq \mathbb E^{\P_0^t}\left[U\big(X_T^{t,x,u}+f\iota\big)-\delta(\iota)-\int_t^T\delta(u_s)ds\right],
\end{equation}
and the same statement holds for any $\iota$ which is bounded and $\Fc^t_{T-}-$measurable by arbitrariness of $\eps_o$. Now, since the map $(x,\iota)\longmapsto U(x+f\iota)-\delta(\iota)$ is Borel measurable, we can argue as in the proof of Proposition 3.1 in \cite{bouchard2014regularity} to obtain the existence for any $\varepsilon>0$ of a Borel measurable map $x\longmapsto \iota_\eps(x)$ such that 
$$\widehat U(X_T^{t,x,u})\leq U\big(X^{t,x,u}_T+f\iota_\eps(X^{t,x,u}_T)\big)-\delta\big(\iota_\eps(X^{t,x,u}_T)\big)+\eps.$$
Then, if we define
$$\iota_{n,\eps}:=\iota_\eps(X^{t,x,u}_T){\bf 1}_{|\iota_\eps(X^{t,x,u}_T)|\leq n},$$
$\iota_{\eps,n}$ is bounded and $\mathcal F_{T-}-$measurable by continuity of the paths $X^{t,x,u}$. Hence by \eqref{eq:face} we have, using the fact that $\delta_T$ is null at $0$
$$Y_t^x\geq \mathbb E^{\P_0^t}\left[\widehat U(X_T^{t,x,u}){\bf 1}_{|\iota_\eps(X^{t,x,u}_T)|\leq n}+U\big(X_T^{t,x,u}\big){\bf 1}_{|\iota_\eps(X^{t,x,u}_T)|> n}-\int_t^T\delta_s(u_s)ds\right]-\eps.$$
Then the required result follows by letting first $n$ go to infinity and dominated convergence (remember that $\widehat U$ is Lipschitz and $X^{t,x,u}$ has moments of any order), and then $\eps$ go to $0$.
\ep

\vspace{0.5em}
\proof[Proof of Lemma \ref{lemma:facelift}]
First of all, one inequality is trivial by taking a constant control $\iota=0$. Then, the following dynamic programming principle holds classically for any $0\leq t\leq s\leq T$
\begin{equation}\label{eq:dpp}
Y_t^{x}=\underset{u\in\mathcal V_{\rm b}^t}{\sup}\mathbb E^{\P_0^t}\left[Y_s^{X^{t,x,u}_s}-\int_t^s\delta(u_r)dr\right].
\end{equation}
Fix now some $(t,x)\in[0,T]\times\mathbb R^d$, some $\iota\in \mathbb L_\infty(\Fc_t)$, some $\eps>0$ small enough, and define
$$u^\eps:=\frac1\eps\iota{\bf 1}_{[t,t+\eps]}\in\mathcal V^t_{\rm b}.$$

By the dynamic programming principle, we have, 
\begin{equation}\label{eq:dyn}
Y_t^{x}\geq \mathbb E^{\P_0^t}\left[Y_{t+\eps}^{X^{t,x,u^\eps}_{t+\eps}}-\int_t^{t+\eps}\delta(u^\eps_r)dr\right]= \mathbb E^{\P_0^t}\left[Y_{t+\eps}^{X_{t+\eps}^{t,x,u^\eps}}-\frac1\eps\int_{t}^{t+\eps}\delta(\iota)dr\right].
\end{equation}
Next, by Lemma \ref{lemma:face} and the tower property we have for any $u\in\mathcal V_{\rm b}^{t}$
$$\mathbb E^{\P_0^t}\left[Y_{t+\eps}^{X_{t+\eps}^{t,x,u^\eps}}\right]\geq \mathbb E^{\mathbb P_0^t}\left[\widehat U\Big(X_T^{t+\eps,X_{t+\eps}^{t,x,u^\eps},u^{t+\eps,B^t}}\Big)-\int_{t+\eps}^T\delta(u_r)dr\right].$$
Be definition of $u^\eps$, we have, using similar arguments as in the proof of Lemma \ref{lemma:face}, that along a subsequence if necessary, 
$$\mathbb E^{\mathbb P_0^t}\left[\widehat U\Big(X_T^{t+\eps,X_{t+\eps}^{t,x,u^\eps},u^{t+\eps,B^t}}\Big)-\int_{t+\eps}^T\delta(u_r)dr\right]\underset{\eps\longrightarrow 0}{\longrightarrow}\mathbb E^{\mathbb P_0^t}\left[\widehat U\Big(X_T^{t,x+f\iota,u}\Big)-\int_{t}^T\delta(u_r)dr\right].$$
Thus, we deduce from passing to the limit in \eqref{eq:dyn} that
$$Y_t^{x}\geq \mathbb E^{\mathbb P_0^t}\left[\widehat U\Big(X_T^{t,x+f\iota,u}\Big)-\delta(\iota)-\int_{t}^T\delta(u_r)dr\right],$$
which implies by Lemma \ref{lemma:face} and arbitrariness of $u\in\mathcal V_{\rm b}^{t}$
$$Y_t^{x}\geq Y_t^{x+f\iota}-\delta(\iota),$$
hence the desired result.
\ep

\vspace{0.5em}
\proof[Proof of Proposition \ref{prop:prop}]
The first result is an immediate consequence of Lemma \ref{lemma:face}, the fact that $\widehat U$ is Lipschitz continuous, and classical estimates on the solutions to SDE satisfied by $X^{t,x,u}$. 

\vspace{0.5em}
Next, we have by \eqref{eq:dpp}, the regularity in $x$ we just proved and \eqref{eq: estimate SDE Delta initial}
\begin{align*}
Y_t^{x}\geq \mathbb E^{\P_0^t}\left[Y_s^{X_s^{t,x,0}}\right]&\geq \mathbb E^{\P_0^t}\left[Y_s^{x}\right]-C\mathbb E^{\P_0^t}\left[\No{X_s^{t,x,0}-x}\right]\geq \mathbb E^{\P_0^t}\left[Y_s^{x}\right]-C\big(1+\No{x}\big)(s-t)^{1/2}.
\end{align*}

Then, we have for any $u\in\mathcal V^t_{\rm b}$, using the fact that $t\longmapsto\delta_t(\cdot)$ is non-increasing and sub--linear, as well as Lemma \ref{lemma:facelift}
\begin{align*}
\mathbb E^{\mathbb P_0^s}\left[U\bigg(X_T^{s,X_s^{t,x,u},u^{s,B^t}}\bigg)-\int_s^T\delta(u_r)dr\right]-\int_t^s\delta(u_r)dr\leq Y_s^{X_s^{t,x,u}}-\int_t^s\delta(u_r)dr&\leq Y_s^{X_s^{t,x,u}}-\delta\left(\int_t^su_rdr\right)\\
&\leq Y_s^{X_s^{t,x, u}-f\int_t^su_rdr}.
\end{align*}
Then by the tower property, we have, taking expectations under $\P_0^t$ on both sides
\begin{align*}
\mathbb E^{\mathbb P_0^t}\left[U(X_T^{t,x,u})-\int_t^T\delta(u_r)dr\right]&\leq \mathbb E^{\P_0^t}\left[Y_s^{X_s^{t,x, u}-f\int_t^su_rdr}\right]\\
&\leq  \mathbb E^{\P_0^t}\left[Y_s^{x}\right]+\mathbb E^{\P_0^t}\left[\No{X_s^{t,x, u}-f\int_t^su_rdr-x}\right]\\
&\leq \mathbb E^{\P_0^t}\left[Y_s^{x}\right]+C\big(1+\No{x}\big)(s-t)^{1/2},
\end{align*}
where we used the fact that in the expression $X_s^{t,x, u}-f\int_t^su_rdr-x$ the control $u$ actually disappears. By definition of $Y_t^x$, this ends the proof.
\ep
\end{document}